\documentclass{article}
\usepackage{amssymb}
\usepackage{amsfonts}
\usepackage{amsmath}

\newtheorem{theorem}{Theorem}

\newtheorem{corollary}[theorem]{Corollary}

\newtheorem{lemma}[theorem]{Lemma}

\begin{document}

\title{Limit theorems for critical branching processes in an extremely
unfavorable random environment\thanks{%
This work was supported by the Russian Science Foundation under grant
no.24-11-00037 https://rscf.ru/en/project/24-11-00037/}}
\author{Vladimir Vatutin\thanks{%
Steklov Mathematical Institute of Russian Academy of Sciences, 8 Gubkina
St., Moscow 119991 Russia Email: vatutin@mi-ras.ru},\, Elena Dyakonova%
\thanks{%
Steklov Mathematical Institute of Russian Academy of Sciences, 8 Gubkina
St., Moscow 119991 Russia Email: elena@mi-ras.ru}}
\date{}
\maketitle

\begin{abstract}
Let $\{Z_{m},m\geq 0\}$ be a critical branching process in random
environment and $\{S_{m},m\geq 0\}$ be its associated random walk. Assuming
that the increments distribution of the associated random walk belongs
without centering to the domain of attraction of an $\alpha $-stable law we
prove conditional limit theorems describing, as $n\rightarrow \infty $, the
distribution the number of particles in the process $\{Z_{m},0\leq m\leq n\}$
given $Z_{n}>0$ and $S_{n}\leq const$.

\textbf{Key words}: stable random walks, branching processes, survival
probability, extreme random environment

\textbf{AMS subject classification}: Primary 60G50, Secondary 60J80, 60K37
\end{abstract}

\section{Introduction and main results}

We study the population size of a critical branching process evolving in an
unfavorable random environment. To describe the problems we plan to consider
denote by $\mathfrak{F}$ $=\left\{ \mathfrak{f}\right\} $ the space of all
probability measures on the set $\mathbb{N}:\mathbb{=}\{1,2,...\}$ and put $%
\mathbb{N}_{0}:=\left\{ 0\right\} \cup \mathbb{N}$. For notational reason,
we identify a measure $\mathfrak{f}=\left\{ \mathfrak{f}(\left\{ 0\right\} ),%
\mathfrak{f}(\left\{ 1\right\} ),...\right\} \in $ $\mathfrak{F}$ with the
respective probability generating function%
\begin{equation*}
f(z)=\sum_{k=0}^{\infty }\mathfrak{f}(\left\{ k\right\} )z^{k},\quad z\in
\lbrack 0,1],
\end{equation*}%
and make no difference between $\mathfrak{f}$ and $f$. Equipped with the
metric of total variation, $\mathfrak{F}$ $=\left\{ \mathfrak{f}\right\}
=\left\{ f\right\} $ becomes a polish space. Let
\begin{equation*}
F(z)=\sum_{j=0}^{\infty }F\left( \left\{ j\right\} \right) z^{j},\quad z\in
\lbrack 0,1],
\end{equation*}%
be a random variable taking values in $\mathfrak{F}$, and let
\begin{equation*}
F_{n}(z)=\sum_{j=0}^{\infty }F_{n}\left( \left\{ j\right\} \right)
z^{j},\quad z\in \lbrack 0,1],\quad n\in \mathbb{N},
\end{equation*}%
be a sequence of independent probabilistic copies of $F$. The infinite
sequence $\mathcal{E}=\left\{ F_{n},n\in \mathbb{N}\right\} $ is called a
random environment.

A sequence of nonnegative random variables $\mathcal{\{}%
Z_{0},Z_{1},Z_{2},...\}$ specified on\ a probability space $(\Omega ,%
\mathcal{F},\mathbf{P})$ is called a branching process in random environment
(BPRE), if $Z_{0}\in \mathbb{N}$ is independent of $\mathcal{E}$ and given $%
\mathcal{E}$ the process $\mathcal{Z}$ is a Markov chain with
\begin{equation*}
\mathcal{L}\left( Z_{n}|Z_{n-1}=k_{n-1},\mathcal{E}=(f_{1},f_{2},...)\right)
=\mathcal{L}(\xi _{n1}+\ldots +\xi _{nk_{n-1}})
\end{equation*}%
for all $n\in \mathbb{N}$, $k_{n-1}\in \mathbb{N}_{0}$ and $%
f_{1},f_{2},...\in \mathfrak{F}$, where $\xi _{n1},\xi _{n2},\ldots $ is a
sequence of i.i.d. random variables with distribution $f_{n}$. Given $%
\mathcal{E}=\left\{ F_{n},n\in \mathbb{N}\right\} $ the sequence
\begin{equation*}
S_{0}=0,\quad S_{n}=X_{1}+...+X_{n},\ n\geq 1,
\end{equation*}%
where $X_{i}=\log F_{i}^{\prime }(1),i=1,2,...,$ is called the associated
random walk for the brancning process in random environment. Under our
assumptions $\left\{ X_{i},i\in \mathbb{N}\right\} $ is a sequence of i.i.d.
random variables.

\bigskip Let
\begin{equation*}
A:=\{\alpha \in (0,2)\backslash \{1\},\,|\beta |<1\}\cup \{\alpha =1,\beta
=0\}\cup \{\alpha =2,\beta =0\}
\end{equation*}%
be a subset in $\mathbb{R}^{2}.$ For $(\alpha ,\beta )\in A$ and a random
variable $X$ we write $X\in \mathcal{D}\left( \alpha ,\beta \right) $ if the
distribution of $X$ belongs to the domain of attraction of a stable law with
density $g_{\alpha ,\beta }(x),x\in (-\infty ,+\infty ),$ and the
characteristic function%
\begin{equation*}
G_{\alpha ,\beta }(w)=\int_{-\infty }^{+\infty }e^{iwx}g_{\alpha .\beta
}(x)\,dx=\exp \left\{ -c|w|^{\,\alpha }\left( 1-i\beta \frac{w}{|w|}\tan
\frac{\pi \alpha }{2}\right) \right\} ,\ c>0,
\end{equation*}%
and, in addition, $\mathbf{E}X=0$ if this moment exists. Note that the law
corresponding to $G_{\alpha ,\beta }(w)$ is not one-sided.

\bigskip We impose of the following restrictions on the properties of the
random environment.

\medskip \textbf{Condition B0.} \emph{The random variables }$X_{n},n\in
\mathbb{N},$\emph{\ are independent copies of a random variable }$X\in
D\left( \alpha ,\beta \right) $\emph{\ whose distribution is under }$\mathbf{%
P}$\emph{\ non-lattice}.

\medskip

Sometimes we will need a stronger assumption:

\medskip

\textbf{Condition B1.} \emph{The random variables }$X_{n},n\in \mathbb{N},$%
\emph{\ are independent copies of a random variable }$X\in D\left( \alpha
,\beta \right) $\emph{\ whose distribution is\ under }$\mathbf{P}$\emph{\
absolutely continuous with respect to the Lebesgue measure on }$\mathbb{R}$%
\emph{, and there exists }$n\in \mathbb{N}$\emph{\ such that the density }$%
f_{n}(x):=\mathbf{P}(S_{n}\in dx)/dx$\emph{\ of }$S_{n}$\emph{\ is bounded.}

\bigskip According to the classification of BPRE's (see, for instance, \cite%
{agkv} and \cite{KV2017}), Condition $B0$ (as well as Condition $B1$) means
that we consider the critical BPRE's.

If \ Condition $B0$ is satisfied, then there is an increasing sequence of
positive numbers
\begin{equation*}
a_{n}\ =\ n^{1/\alpha }\ell (n)
\end{equation*}%
with a slowly varying sequence $\ell (1),\ell (2),\ldots ,$ such that
\begin{equation*}
\left\{ \frac{S_{\left[ nt\right] }}{a_{n}},t\geq 0\right\} \Longrightarrow
\mathcal{Y}=\left\{ Y_{t},t\geq 0\right\}
\end{equation*}%
as $n\rightarrow \infty ,$ where%
\begin{equation*}
\mathbf{E}e^{iwY_{t}}=G_{\alpha ,\beta }(wt^{1/\alpha }),\quad t\geq 0,
\end{equation*}%
and the symbol $\Longrightarrow $ stands for the weak convergence in the
space $D[0,\infty )$ of c\`{a}dl\`{a}g functions endowed with Skorokhod
topology. In addition,
\begin{equation*}
\lim_{n\rightarrow \infty }\mathbf{P}\left( S_{n}>0\right) =\rho =\mathbf{P}%
\left( Y_{1}>0\right) \in (0,1).
\end{equation*}

\bigskip Our third assumption on the environment concerns reproduction laws
of particles. Set%
\begin{equation*}
\gamma (b):=\frac{\sum_{k=b}^{\infty }k^{2}F\left( \left\{ k\right\} \right)
}{\left( \sum_{i=0}^{\infty }iF\left( \left\{ i\right\} \right) \right) ^{2}}%
.
\end{equation*}

\paragraph{Condition B2.}

\emph{There exist $\varepsilon >0$ and $b\in \,$}$\mathbb{N}$ \emph{such
that\ } \emph{\ }
\begin{equation*}
\mathbf{E}[(\log ^{+}\gamma (b))^{\alpha +\varepsilon }]\ <\ \infty ,
\end{equation*}%
\emph{where }$\log ^{+}x=\log (x\vee 1)$\emph{.}

In what follows we assume, if otherwise\ is not stated that $Z_{0}=1$. It is
known (see \cite[Theorem 1.1 and Corollary 1.2]{agkv}) that if Conditions $B0
$ and $B2$ are valid, then there exist a number $\theta \in (0,\infty )$ and
a sequence $l(1),l(2),...,$ slowly varying at infinity such that, as $%
n\rightarrow \infty $%
\begin{equation*}
\mathbf{P}\left( Z_{n}>0\right) \sim \theta \mathbf{P}\left( L_{n}\geq
0\right) \sim \theta n^{-(1-\rho )}l(n),
\end{equation*}%
and, for any $x\geq 0$
\begin{eqnarray}
\mathbf{P}\left( Z_{n}>0,S_{n}\leq xa_{n}\right)  &=&\mathbf{P}\left(
S_{n}\leq xa_{n}|Z_{n}>0\right) \mathbf{P}\left( Z_{n}>0\right)   \notag \\
&\sim &\mathbf{P}\left( Y_{1}^{+}\leq x\right) \mathbf{P}\left(
Z_{n}>0\right) ,  \label{Meander00}
\end{eqnarray}%
where $\mathcal{Y}^{+}=\left\{ Y_{t}^{+},0\leq t\leq 1\right\} $ denotes the
meander of a strictly stable process $\mathcal{Y}$ with index $\alpha $,
i.e. the strictly $\alpha $-stable Levy process staying positive on the
seminterval $(0,1]$ (see \cite{Do85}, \cite{Du78}).

Since $\mathbf{P}\left( Y_{1}^{+}\leq 0\right) =0,$ it follows from (\ref%
{Meander00}) that if a positive function $\varphi (n),n\in \mathbb{N},$
satisfies the restriction
\begin{equation*}
\limsup_{n\rightarrow \infty }\frac{\varphi (n)}{a_{n}}\leq 0,
\end{equation*}%
then%
\begin{equation*}
\mathbf{P}\left( Z_{n}>0,S_{n}\leq \varphi (n)\right) =\mathbf{P}\left(
S_{n}\leq \varphi (n)|Z_{n}>0\right) \mathbf{P}\left( Z_{n}>0\right)
=o\left( \mathbf{P}\left( Z_{n}>0\right) \right)
\end{equation*}%
as $n\rightarrow \infty $.

For this reason we consider the environment meeting the condition $%
S_{n}=o(a_{n})$ as unfavorable for the development of the critical BPRE's.

An important case of the unfavorable random environment was considered in
\cite{VD2022} and \cite{VDD2024}. To formulate some results of the papers we
set%
\begin{eqnarray}
M_{n} &:&=\max \left( S_{1},...,S_{n}\right) ,\quad L_{n}:=\min \left(
S_{1},...,S_{n}\right) ,  \notag \\
\tau _{n} &:&=\min \left\{ 0\leq k\leq n:S_{k}=\min (0,L_{n})\right\} ,
\notag
\end{eqnarray}%
\begin{equation}
b_{n}:=\frac{1}{a_{n}n}=\frac{1}{\ n^{1/\alpha +1}\ell (n)},  \label{Defb}
\end{equation}%
and introduce two renewal functions
\begin{eqnarray*}
U(x):= &&\mathbf{I}\left\{ x\geq 0\right\} +\sum_{n=1}^{\infty }\mathbf{P}%
\left( S_{n}\geq -x,M_{n}<0\right)  \\
&=&\mathbf{I}\left\{ x\geq 0\right\} +\sum_{n=1}^{\infty }\mathbf{P}\left(
S_{n}\geq -x,\tau _{n}=n\right)
\end{eqnarray*}%
and
\begin{equation*}
V(x):=\mathbf{I}\left\{ x<0\right\} +\sum_{n=1}^{\infty }\mathbf{P}\left(
S_{n}<-x,L_{n}\geq 0\right) ,
\end{equation*}%
where $\mathbf{I}\left\{ A\right\} $ is the indicator of the event $A$.

It was shown in \cite{VD2022} that if Conditions $B1$ and $B2$ are valid and
$\varphi (n)\rightarrow \infty $ as $n\rightarrow \infty $ in such a way
that $\varphi (n)=o(a_{n}),$ then
\begin{equation}
\mathbf{P}\left( Z_{n}>0,S_{n}\leq \varphi (n)\right) \sim Dg_{\alpha ,\beta
}(0)b_{n}\int_{0}^{\varphi (n)}V(-u)du,\;D\in (0,\infty ).  \label{Lim_surv0}
\end{equation}%
The same assumptions were used in \cite{VD2023} and \cite{VDD2023} to find
the limits, as $n\rightarrow \infty ,$ of the conditional distributions of
the properly scaled processes $\left\{ Z_{m},0\leq m\leq n\right\} $ and $%
\left\{ \log Z_{m},0\leq m\leq n\right\} $ given the event $\left\{
Z_{n}>0,S_{n}\leq \varphi \left( n\right) \right\} $.

These results were complemented in \cite{VDD2024} by proving that if $%
S_{n}\leq K$ for a fixed constant $K$ (we call such an environment \textit{%
extremely unfavorable}) and Conditions $B1$ and $B2$ are met, then

\begin{equation}
\mathbf{P}\left( Z_{n}>0,S_{n}\leq K\right) \sim b_{n}\left(
G_{left}(K)+G_{right}(K)\right)  \label{Limit_surviv}
\end{equation}%
as $n\rightarrow \infty ,$ where $G_{left}(K)$ and $G_{right}(K)$ are
positive and finite constants whose explicit forms are given by formulas (%
\ref{Gleft0}) and (\ref{Gright0}) below.

In this paper we use the results of \cite{VDD2024} to describe, as $%
n\rightarrow \infty $ properties of the trajectories of the process $\left\{
Z_{m},0\leq m\leq n\right\} $ given $\left\{ Z_{n}>0,S_{n}\leq K\right\} $.

Our first main result concerns the conditional distribution of $Z_{n}$.

\begin{theorem}
\label{T_simplelaw} If conditions $B1$ and $B2$ are valid, then the
conditional laws $\mathcal{L}(Z_{n}|Z_{n}>0,S_{n}\leq K),\,n\geq 1,$
converge weakly, as $n\rightarrow \infty $ to a proper discrete probability
distribution.
\end{theorem}

The second main statement of the present paper describes the limiting
behavior of the rescaled generation size process. To formulate the claimed
theorem we introduce, for integers $m\geq k$ the notation%
\begin{equation}
\mathcal{Y}^{k,m}=\{Y^{k,m}(t),\quad t\in \lbrack 0,1]\},  \label{Def_Y}
\end{equation}%
where%
\begin{equation*}
Y^{k,m}(t):=\exp \left\{ -S_{k+\lfloor \left( m-2k\right) t\rfloor }\right\}
Z_{k+\lfloor (m-2k)t\rfloor }.
\end{equation*}

Thus, $\mathcal{Y}^{0,m}$ is a tragectory of the rescaled population size of
the branching process on the time-interval $[0,m]$ and $\mathcal{Y}^{k,m}$
is a piece of the tragectory between the time-points $k$ and $m-k$.

\begin{theorem}
\label{T_generalcase} If Conditions $B1$ and $B2$ are valid then, for any $%
\theta \in \left( 0,1/2\right) $
\begin{equation*}
\mathcal{L}(\mathcal{Y}^{\theta n,n}|Z_{n}>0,S_{n}\leq K)\Rightarrow
\mathcal{L}(W(t),t\in \lbrack 0,1])
\end{equation*}%
as $n\rightarrow \infty $, where the limiting process $\mathcal{W}%
=\{W(t),t\in \lbrack 0,1]\}$ is a stochastic process with a.s. continuous
paths, that is $\mathbf{P}\left( W(t)=W\text{ for all }t\in \lbrack
0,1]\right) =1$ for some random variable $W$. Furthermore,
\end{theorem}

\begin{equation*}
\mathbf{P}\left( 0<W<\infty \right) =1.
\end{equation*}

Here and in what follows we consider expressions of the form $\theta
n,(1-2\theta )n$ as $\left[ \theta n\right] ,\left[ (1-2\theta )n\right] $
and use the symbol $\Rightarrow $ to denote weak convergence of c\`{a}dl\`{a}%
g functions with respect to Skorokhod topology in the space $D[0,1]$.

The structure of the remaining part of the paper looks as follows. In
Section \ref{Sec2} we formulate a number of known results for random walks
conditioned to stay nonnegative or negative. Section \ref{sec3} is devoted
to the proof of Theorem \ref{T_simplelaw}. In Section \ref{sec4} we proof
Theorem \ref{T_generalcase}.

In what follows we denote by $C,C_{1},C_{2},...,$ some positive constants
that may be different in different formulas or even within one and the same
formula.

\section{Auxiliary results\label{Sec2}}

By means of the renewal functions $U$ and $V$ we introduce two new
probability measures $\mathbf{P}^{+}$ and $\mathbf{P}^{-}$ using the
identities
\begin{equation}
\begin{array}{rl}
\mathbf{E}[U(x+X);X+x\geq 0]\ =\ U(x), & x\geq 0, \\
\mathbf{E}[V(x+X);X+x<0]\ =\ V(x), & x\leq 0,%
\end{array}
\label{Martingale}
\end{equation}%
valid for any oscillating random walk (see \cite{KV2017}, Chapter 4.4.3).
The construction procedure of these measures is explained in detail in \cite%
{agkv} and \cite{ABGV2011} (see also \cite{KV2017}, Chapter 5.2). We recall
here only some basic definitions related to this construction.

Let $\mathcal{F}_{n}$ be the $\sigma $-algebra of events generated by random
variables $F_{1},F_{2},...,F_{n}$ and $Z_{0},Z_{1},...,Z_{n}$. The sequence $%
\left\{ \mathcal{F}_{n},n\geq 1\right\} $ forms a filtration $\mathcal{F}$.
We assume that the random walk $\mathcal{S=}\left\{ S_{n},n\geq 0\right\} $
with the initial value $S_{0}=x,\,x\in \mathbb{R}$, is adapted to the
filtration $\mathcal{F}$ and construct for $x\geq 0$ probability measure $%
\mathbf{P}_{x}^{+}$ and expectation $\mathbf{E}_{x}^{+}$ as follows. For
every sequence $T_{0},T_{1},...,$ of random variables with values in some
space $\mathcal{T}$ and adopted to $\mathcal{F}$, and for any bounded and
measurable function $g:\mathcal{T}^{n+1}\rightarrow \mathbb{R}$, $n\in
\mathbb{N}_{0}$, we set%
\begin{equation*}
\mathbf{E}_{x}^{+}[g(T_{0},\ldots ,T_{n})]\ :=\ \frac{1}{U(x)}\mathbf{E}%
_{x}[g(T_{0},\ldots ,T_{n})U(S_{n});L_{n}\geq 0].\
\end{equation*}%
Similarly, for $x\leq 0$ we set
\begin{equation*}
\mathbf{E}_{x}^{-}[g(T_{0},\ldots ,T_{n})]\ :=\ \frac{1}{V(x)}\mathbf{E}%
_{x}[g(T_{0},\ldots ,T_{n})V(S_{n});M_{n}<0]\ .
\end{equation*}

In view of relations (\ref{Martingale}) the definitions above are correct
and adopted with respect to $n$. We write, for brevity, $\mathbf{P}^{\pm }$
and $\mathbf{E}^{\pm }$ for $\mathbf{P}_{0}^{\pm }$ and $\mathbf{E}_{0}^{\pm
}$.

We now formulate, for reader's convenience, several known lemmas playing an
important role in proving Theorems \ref{T_simplelaw} and \ref{T_generalcase}.

The first two lemmas are related to the properties of random walks only.

\begin{lemma}
\label{L_asymptot} \emph{If Condition }$\emph{B0}$ is valid, then (see \cite%
{VW2010}, Theorem 4)
\begin{equation}
\mathbf{P}\left( S_{n}\leq y,L_{n}\geq 0\right) \sim g_{\alpha ,\beta
}(0)b_{n}\int_{0}^{y}V(-u)du\;  \label{AsymS_tay}
\end{equation}%
as $n\rightarrow \infty $ uniformly in $y\in (0,\delta _{n}a_{n}]$, where $%
\delta _{n}\rightarrow 0$ as $n\rightarrow \infty $. Besides (see \cite%
{VD2022}, Lemma 2), there is a constant $C\in (0,\infty )$ such that%
\begin{equation}
\mathbf{P}\left( S_{n}\leq y,L_{n}\geq 0\right) \leq
Cb_{n}\int_{0}^{y}V(-u)du  \label{EstimS_tay}
\end{equation}%
for all $y>0$ and all $n\in \mathbb{N}$.
\end{lemma}

It follows from (\ref{Lim_surv0}) and (\ref{AsymS_tay}) that%
\begin{equation*}
\mathbf{P}\left( Z_{n}>0,S_{n}\leq y\right) \sim D\mathbf{P}\left( S_{n}\leq
y,L_{n}\geq 0\right)
\end{equation*}%
as $n\rightarrow \infty $ uniformly in $y\in (0,\delta _{n}a_{n}]$ , where $%
\delta _{n}\rightarrow 0$ as $n\rightarrow \infty $.

\begin{lemma}
\label{L_cond} (see Lemma 4 in \cite{VD2022}) Assume Condition $B0$. Let $%
H_{1},H_{2},...,$ be a uniformly bounded sequence of real-valued random
variables adopted to some filtration $\mathcal{\tilde{F}=}\left\{ \mathcal{%
\tilde{F}}_{k},k\in \mathbb{N}\right\} $, which converges $\mathbf{P}^{+}$%
-a.s. to a random variable $H_{\infty }$ as $n\rightarrow \infty $. Suppose
that $\varphi (n),$ $n\in \mathbb{N},$ is a real-valued function such that $%
\inf_{n\in \mathbb{N}}\varphi (n)\geq C>~0$ and $\varphi (n)=o(a_{n})$ as $%
n\rightarrow \infty $. Then%
\begin{equation*}
\lim_{n\rightarrow \infty }\frac{\mathbf{E}\left[ H_{n};S_{n}\leq \varphi
(n),L_{n}\geq 0\right] }{\mathbf{P}\left( S_{n}\leq \varphi (n),L_{n}\geq
0\right) }=\mathbf{E}^{+}\left[ H_{\infty }\right] .
\end{equation*}
\end{lemma}

The next statement is borrowed from \cite{agkv}. We here write it in a more
detailed way.

\begin{lemma}
\label{Lemprop3.1} (see Proposition 3.1 in \cite{agkv}) \emph{Let }%
Conditions $B0$ and $B2$ be valid. Then for any $q=1,2,...$
\begin{equation*}
\mathbf{P}^{+}\left( Z_{n}>0,Z_{0}>q\text{ for all }n|\mathcal{E}\right)
>0\quad \;\mathbf{P}^{+}-\text{a.s.}
\end{equation*}%
In particular,%
\begin{equation*}
\mathbf{P}^{+}\left( Z_{n}>0,Z_{0}>q\text{ for all }n\right) >0.
\end{equation*}%
Moreover,%
\begin{equation}
e^{-S_{n}}Z_{n}\mathbf{I}\left\{ Z_{0}=q\right\} \rightarrow W_{\left(
q\right) }^{+}\;\quad \mathbf{P}^{+}-\text{a.s.}  \label{Funct}
\end{equation}%
as $n\rightarrow \infty $, where the random variable $W_{\left( q\right)
}^{+}$ is such that%
\begin{equation}
\left\{ W_{\left( q\right) }^{+}>0\right\} =\left\{ Z_{0}=q,Z_{n}>0\text{
for all }n\right\} \quad \;\mathbf{P}^{+}-\text{a.s.}  \label{DefWPlus}
\end{equation}
\end{lemma}

\bigskip For the sake of brevity we write in the sequel $W^{+}$ for $%
W_{\left( 1\right) }^{+}.$

\begin{lemma}
\label{L_Subcr} (see Lemma 5 in \cite{VDD2024}) \emph{Let }Conditions $B0$
and $B2$ be valid and $0<\delta <1$ be fixed. Let
\begin{equation*}
A_{n}=\mathfrak{a}_{n}(F_{1},\ldots ,F_{\lfloor \delta n\rfloor
},Z_{0},Z_{1},...,Z_{\lfloor \delta n\rfloor }),\ n\geq 1,
\end{equation*}%
be random variables with values in an Euclidean (or polish) space $\mathcal{A%
}$ such that, for all $x\geq 0$%
\begin{equation*}
A_{n}\ \rightarrow A_{\infty }\quad \mathbf{P}_{x}^{+}\text{-a.s.}
\end{equation*}%
for some $\mathcal{A}$-valued random variable $A_{\infty }$. Also let $B_{n}=%
\mathfrak{b}_{n}(F_{1},\ldots ,F_{\lfloor \delta n\rfloor }),n\geq 1$, be
random variables with values in an Euclidean (or polish) space $\mathcal{B}$
such that
\begin{equation*}
B_{n}\ \rightarrow \ B_{\infty }\quad \mathbf{P}^{-}\text{-a.s.}
\end{equation*}%
for some $\mathcal{B}$-valued random variable $B_{\infty }.$ Denote
\begin{equation*}
\tilde{B}_{n}\ :=\ \mathfrak{b}_{n}(F_{n},\ldots ,F_{n-\lfloor \delta
n\rfloor +1})\ .
\end{equation*}%
Then for any bounded, continuous function $\vec{\Psi}:\mathcal{A}\times
\mathcal{B}\times \mathbb{R}\rightarrow \mathbb{R}$ and any $\eta >0$
\begin{align*}
& \frac{\mathbf{E}[\vec{\Psi}(A_{n},\tilde{B}_{n},S_{n})e^{\eta
S_{n}};\;\tau _{n}=n]}{\mathbf{E}[e^{\eta S_{n}};\tau _{n}=n]}\; \\
& \qquad \qquad \qquad \rightarrow \ \iiint \vec{\Psi}(\mathfrak{a},%
\mathfrak{b},-y)\,\mathbf{P}_{y}^{+}\left( A_{\infty }\in d\mathfrak{a}%
\right) \mathbf{P}^{-}\left( B_{\infty }\in d\mathfrak{b}\right) \mu _{\eta
}(dy)
\end{align*}%
as $n\rightarrow \infty $, where%
\begin{equation}
\mu _{\eta }(dy)=\frac{e^{-\eta y}U(y)\,}{\int_{0}^{\infty }e^{-\eta
x}U(x)\,dx}dy,\ \eta >0.  \label{Def_mu}
\end{equation}
\end{lemma}

In the sequel we apply Lemma \ref{L_Subcr} to some spaces related to the
space $D\mathcal{[}0,1\mathcal{]}$ of c\`{a}dl\`{a}g functions. For this
reason we recall some basic fact concerning this space. \ Denote by $\Lambda
$ the set of all strictly increasing and continuous function $\lambda
(t),t\in \lbrack 0,1],$ with $\lambda (0)=0$ and $\lambda (1)=1$. For $\phi
_{1}(t)\in D\mathcal{[}0,1\mathcal{]}$ and $\phi _{2}(t)\in D\mathcal{[}0,1%
\mathcal{]}$ we define two metrics:

Skorohod $d$ metric:%
\begin{equation*}
d(\phi _{1},\phi _{2})=\inf_{\lambda \in \Lambda }\max \left\{ \sup_{t\in
\lbrack 0,1]}\left\vert \lambda (t)-t\right\vert ,\sup_{t\in \lbrack
0,1]}\left\vert \phi _{1}(\lambda (t))-\phi _{2}(t)\right\vert \right\}
\end{equation*}%
and $d^{0}$ metric:%
\begin{equation*}
d^{0}(\phi _{1},\phi _{2})=\inf_{\lambda \in \Lambda }\max \left\{
\sup_{s<t}\log \frac{\lambda (t)-\lambda (s)}{t-s},\sup_{t\in \lbrack
0,1]}\left\vert \phi _{1}(\lambda (t))-\phi _{2}(t)\right\vert \right\} .
\end{equation*}%
It is known (see, for instance, \cite{Bil99}, Section 12, Theorem 12.1) that
metrics $d$ and $d^{0}$ are equivalent and, moreover, $D\mathcal{[}0,1%
\mathcal{]}$ is separable and complete under $d^{0}$, and, therefore,
polish. These results lead to the following statement.

\begin{corollary}
\label{C_polish} The space $\mathcal{A}:\mathcal{=}D\mathcal{[}0,1\mathcal{%
]\times }\mathbb{R}_{+}$ is polish under the metric $d^{0}\times \left\Vert
x\right\Vert _{\infty }$.
\end{corollary}

For $0\leq z<1$ define iterations
\begin{equation*}
F_{k,n}(z):=F_{k+1}(F_{k+2}(\ldots F_{n}(z)\ldots ))\text{ if }0\leq k<n%
\text{ and }F_{n,n}(z):=z.
\end{equation*}%
The next statement is important in proving Theorem \ref{T_simplelaw}.

\begin{lemma}
\label{L_hbound} Let Condition $B0$ be valid. Then there exists a constant $%
C $ such that%
\begin{equation*}
\mathbf{E}\left[ 1-F_{0,j}(z);S_{j}\leq w,\tau _{j}=j\right] \leq
C(1-z)b_{j}e^{w/2}
\end{equation*}%
for all $j\in \mathbb{N}$ and $w\in \mathbb{R}$.
\end{lemma}

\textbf{Proof}. Since $\left\{ S_{j}\leq w,\tau _{j}=j\right\} =\left\{
S_{j}\leq 0,\tau _{j}=j\right\} $, it is sufficient to consider $w\leq 0$
only. By the duality principle for random walks we obtain
\begin{eqnarray*}
\mathbf{E}\left[ 1-F_{0,j}(z);S_{j}\leq w,\tau _{j}=j\right] &\leq &(1-z)%
\mathbf{E}\left[ e^{S_{j}};S_{j}\leq w,\tau _{j}=j\right] \\
&=&(1-z)\mathbf{E}\left[ e^{S_{j}};S_{j}\leq w,M_{j}<0\right] .
\end{eqnarray*}%
It is shown in \cite[Proposition 2.3]{ABGV2011} that, given Condition $B0$
there exists a constant $C>0$ such that, for all $j\in \mathbb{N}$ and $%
y\leq 0$
\begin{equation*}
\mathbf{P}\left( y\leq S_{j}<y+1,M_{j}<0\right) \leq Cb_{j}U\left( -y\right)
.
\end{equation*}%
Hence it follows that%
\begin{eqnarray*}
\mathbf{E}\left[ e^{S_{j}};S_{j}\leq w,M_{j}<0\right] &\leq &\sum_{i\geq
-w}e^{-i+1}\mathbf{P}\left( -i\leq S_{j}<-i+1,M_{j}<0\right) \\
&\leq &Cb_{j}\sum_{i\geq -w}e^{-i+1}U\left( i\right) \\
&\leq &Cb_{j}e^{w/2+1}\sum_{i\geq -w}e^{-i/2}U\left( i\right) \leq
C_{1}b_{j}e^{w/2},
\end{eqnarray*}%
since%
\begin{equation*}
\int_{0}^{\infty }e^{-\eta x}U(x)dx<\infty
\end{equation*}%
for any $\eta >0$ in view of the inequality $U(x)<C_{2}\left( x+1\right) ,$ $%
x\geq 0,$ valid, for each own $C_{2},$ for any renewal function.

Lemma \ref{L_hbound} is proved.

For $z\in \lbrack 0,1],w\in \mathbb{R},$ and $j\in \mathbb{N}_{0}$ introduce
the notation%
\begin{equation}
O_{j}\left( z,w\right) :=\frac{\mathbf{E}\left[ 1-F_{0,j}(z);S_{j}\leq
w,\tau _{j}=j\right] }{\mathbf{E}\left[ e^{S_{j}};\tau _{j}=j\right] }.
\label{Def_O}
\end{equation}

\begin{corollary}
\label{C_h}\emph{If }Condition $B0$ is valid, then there exists a constant $%
C $ such that%
\begin{equation*}
O_{j}\left( z,w\right) \leq C\left( 1-z\right) e^{w/2}
\end{equation*}%
for all $j\in \mathbb{N}$ and $w\in \mathbb{R}$.
\end{corollary}

\textbf{Proof}. It is known (see, for instance, formula (1.6) in \cite%
{VDD2024} or, in a slightly different notation, Proposition 2.1 in \cite%
{ABGV2011}) that
\begin{equation}
\mathbf{E}\left[ e^{S_{j}};\tau _{j}=j\right] =\mathbf{E}\left[
e^{S_{j}};M_{j}<0\right] \sim g_{\alpha ,\beta }(0)b_{j}\int_{0}^{\infty
}e^{-y}U(y)dy  \label{AsymExponent}
\end{equation}%
as $j\rightarrow \infty $. Using this equivalence and Lemma \ref{L_hbound}
we arrive to the statement of the corollary.

\section{\textbf{Proof of Theorem \protect\ref{T_simplelaw}}\label{sec3}}

We prove Theorem \ref{T_simplelaw} by the arguments similar to those used in
\cite{VDD2024} to investigate the asymptotic behavior of the probability of
the event
\begin{equation*}
\mathcal{R}(n,K):=\left\{ Z_{n}>0,S_{n}\leq K\right\} .
\end{equation*}

We need to show that
\begin{equation*}
\lim_{n\rightarrow \infty }\mathbf{E[}z^{Z_{n}}|\mathcal{R}%
(n,K)]=\lim_{n\rightarrow \infty }\frac{\mathbf{E[}z^{Z_{n}};\mathcal{R}%
(n,K)]}{\mathbf{P}\left( \mathcal{R}(n,K)\right) }
\end{equation*}%
exists for any $z\in \left[ 0,1\right] $ and the limiting distribution has
the properties claimed in the statement of the theorem. To perform this task
we write the decomposition%
\begin{equation*}
\mathbf{E[}z^{Z_{n}};\mathcal{R}(n,K)]=\Lambda \left( 0,J\right) +\Lambda
\left( J+1,n-J\right) +\Lambda \left( n-J+1,n\right) ,
\end{equation*}%
where%
\begin{equation*}
\Lambda \left( A,B\right) :=\mathbf{E[}z^{Z_{n}};\mathcal{R}(n,K),\tau
_{n}\in \lbrack A,B]].
\end{equation*}%
First we note that%
\begin{equation*}
\Lambda \left( J+1,n-J\right) \leq \mathbf{P}\left( \mathcal{R}(n,K),\tau
_{n}\in \left[ J+1,n-J\right] \right)
\end{equation*}%
and
\begin{eqnarray}
&&\limsup_{J\rightarrow \infty }\limsup_{n\rightarrow \infty }\frac{\mathbf{P%
}\left( \mathcal{R}(n,K),\tau _{n}\in \left[ J+1,n-J\right] \right) }{b_{n}}
\notag \\
&&\quad \leq \limsup_{J\rightarrow \infty }\limsup_{n\rightarrow \infty }%
\frac{\mathbf{E}\left[ e^{S_{\tau _{n}}};S_{n}\leq K,\tau _{n}\in \left[
J+1,n-J\right] \right] }{b_{n}}=0  \label{441}
\end{eqnarray}%
according to estimate (7.8) in \cite{VDD2024}. Therefore,

\begin{equation}
\limsup_{J\rightarrow \infty }\limsup_{n\rightarrow \infty }\frac{\Lambda
\left( J+1,n-J\right) }{b_{n}}=0.\quad  \label{IntermSmall}
\end{equation}%
Taking (\ref{Limit_surviv}) into account we see that it remains to only
analyse the case%
\begin{equation*}
\tau _{n}\in \lbrack 0,J]\cup \lbrack n-J+1,n]
\end{equation*}%
for sufficiently large but fixed $J$ and $n\rightarrow \infty $.

We now fix sufficiently large positive integers $N>\left\vert K\right\vert $
and $J$ and select $j\in \lbrack 0,J]$.

It was shown in the proof of Theorem 5 in \cite{VDD2024} that, for any $%
\varepsilon >0$ and all $n\geq j$ there exists $N_{0}=N_{0}(\varepsilon )$
such that
\begin{equation*}
\mathbf{P}\left( \mathcal{R}(n,K),\tau _{n}=j,S_{j}<-N\right) \leq
\varepsilon b_{n-j}
\end{equation*}%
for all $N\geq N_{0}$, and there exist $Q_{0}=Q_{0}(\varepsilon ,N)$ and a
constant $C_{1}$ such that
\begin{equation*}
\mathbf{P}\left( \mathcal{R}(n,K),\tau _{n}=j,S_{j}\geq -N,Z_{j}>Q\right)
\leq C_{1}\varepsilon b_{n-j}
\end{equation*}%
for all $Q\geq Q_{0}$.

Hence it follows that, for any $\varepsilon >0$ and all $n\geq j$
\begin{eqnarray}
&&\mathbf{E}\left[ z^{Z_{n}};\mathcal{R}(n,K),\tau _{n}=j,S_{j}<-N\right]
\notag \\
&&\qquad \leq \mathbf{P}\left( \mathcal{R}(n,K),\tau _{n}=j,S_{j}<-N\right)
\leq \varepsilon b_{n-j}  \label{left1}
\end{eqnarray}%
and
\begin{eqnarray}
&&\mathbf{E}\left[ z^{Z_{n}};\mathcal{R}(n,K),\tau _{n}=j,S_{j}\geq
-N,Z_{j}>Q\right]  \notag \\
&&\qquad \leq \mathbf{P}\left( \mathcal{R}(n,K),\tau _{n}=j,S_{j}\geq
-N,Z_{j}>Q\right) \leq \varepsilon C_{1}b_{n-j}  \label{left2}
\end{eqnarray}%
for sufficiently large $N$ and $Q=Q(\varepsilon ,N)$.

Thus, for any $\varepsilon >0$ and $j\in \lbrack 0,J]$ we have the
representation%
\begin{eqnarray*}
\mathbf{E}\left[ z^{Z_{n}};\mathcal{R}(n,K),\tau _{n}=j\right]
&=&\varepsilon _{j,N,Q}b_{n-j} \\
&+&\mathbf{E}\left[ z^{Z_{n}};\mathcal{R}(n,K),\tau _{n}=j,S_{j}\geq
-N,Z_{j}\in \lbrack 1,Q]\right] ,
\end{eqnarray*}%
where%
\begin{equation*}
\limsup_{N\rightarrow \infty }\limsup_{Q\rightarrow \infty }\sup_{0\leq
j\leq J}\varepsilon _{j,N,Q}=0.
\end{equation*}

We now consider the term%
\begin{equation*}
\mathbf{E}\left[ z^{Z_{n}};\mathcal{R}(n,K),\tau _{n}=j,S_{j}\geq
-N,Z_{j}\in \lbrack 1,Q]\right]
\end{equation*}%
which requires more delicate estimates.

First observe that
\begin{equation*}
\mathbf{E[}z^{Z_{n}}\mathbf{I}\left\{ Z_{n}>0,Z_{0}=q\right\} |\mathcal{E}%
]=F_{0,n}^{q}\left( z\right) -F_{0,n}^{q}\left( 0\right) ,\;z\in (0,1).
\end{equation*}

Further, it is known (see \cite[Theorem 5]{AK1971}) that, for any sequence $%
f_{1},f_{2},...,f_{n},...$ of probability generating functions \
\begin{equation*}
\lim_{n\rightarrow \infty }f_{1}(f_{2}(\ldots f_{n}(z)\ldots
))=\lim_{n\rightarrow \infty }f_{0,n}(z)=:f_{0,\infty }(z)
\end{equation*}%
exists for all $0\leq z<1$. Therefore, for any $q\in \mathbb{N}$
\begin{eqnarray}
&&\lim_{n\rightarrow \infty }\mathbf{E[}z^{Z_{n}}\mathbf{I}\left\{
Z_{n}>0,Z_{0}=q\right\} |\mathcal{E}]=F_{0,\infty }^{q}(z)-F_{0,\infty
}^{q}(0)  \notag \\
&&\qquad \qquad =:\mathbf{E[}z^{Z_{\infty }}\mathbf{I}\left\{ Z_{\infty
}\geq 1,Z_{0}=q\right\} |\mathcal{E}]=:\mathcal{M}(z,q)  \label{Dob501}
\end{eqnarray}%
exists $\mathbf{P}^{+}$-a.s.

It follows from (\ref{Dob501}) and Lemma \ref{L_cond} that, for any $K>0$
\begin{eqnarray*}
&&\mathbf{E}\left[ z^{Z_{n}}\mathbf{I}\left\{ Z_{n}>0,Z_{0}=q\right\}
|S_{n}\leq K,L_{n}\geq 0\right] \rightarrow \mathbf{E}^{+}\left[ \mathcal{M}%
(z,q)\right] \\
&&\qquad \qquad \qquad =\mathbf{E}^{+}\mathbf{[}z^{Z_{\infty }}\mathbf{I}%
\left\{ Z_{\infty }\geq 1,Z_{0}=q\right\} ]=:\mathbf{E}^{+}(z,q)
\end{eqnarray*}%
as $n\rightarrow \infty $. Since $\mathbf{E}\left[ Z_{n}\mathbf{I}\left\{
Z_{0}=q\right\} |\mathcal{E}\right] =qe^{S_{n}}$, we have
\begin{eqnarray*}
&&\mathbf{E}\left[ z^{Z_{n}}\mathbf{I}\left\{ Z_{n}>T,Z_{0}=q\right\}
;S_{n}\leq K,L_{n}\geq 0\right] \\
&&\qquad\qquad\leq\mathbf{P}\left( Z_{n}>T,Z_{0}=q,S_{n}\leq K,L_{n}\geq
0\right) \\
&&\qquad\qquad\leq\frac{q}{T}\mathbf{E}\left[ e^{S_{n}};S_{n}\leq
K,L_{n}\geq 0\right] \\
&&\qquad\qquad\leq\frac{qe^{K}}{T}\mathbf{P}\left( S_{n}\leq K,L_{n}\geq
0\right)
\end{eqnarray*}%
for any $T>0$. Hence we deduce that $Z_{n}\mathbf{I}\left\{ Z_{0}=q\right\} $
does not escape to infinity as $n\rightarrow \infty $ under the condition $%
\left\{ S_{n}\leq K,L_{n}\geq 0\right\} $. Thus,
\begin{eqnarray*}
\lim_{z\uparrow 1}\mathbf{E}^{+}(z,q) &=&\mathbf{E}^{+}(1,q)=\mathbf{P}%
^{+}\left( 1\leq Z_{\infty }<\infty ,Z_{0}=q\right) \\
&=&\mathbf{P}^{+}\left( 1\leq Z_{\infty },Z_{0}=q\right) >0
\end{eqnarray*}%
in view of Lemma \ref{Lemprop3.1} and, therefore, $\mathbf{E}^{+}(z,q)>0$
for any $z\in (0,1]$.

Using now estimate (\ref{AsymS_tay}) we see that, for fixed $q\in \lbrack
1,Q]$ and $x\in \lbrack -N,K]$%
\begin{eqnarray*}
&&\mathbf{E[}z^{Z_{n-j}};Z_{n-j}>0,S_{n-j}\leq K-x,L_{n-j}\geq 0,Z_{0}=q] \\
&&\qquad =\mathbf{E}\left[ F_{0,n-j}^{q}(z)-F_{0,n-j}^{q}(0);S_{n-j}\leq
K-x,L_{n-j}\geq 0\right] \\
&&\qquad \sim \mathbf{E}^{+}(z,q)\mathbf{P}\left( S_{n-j}\leq
K-x,L_{n-j}\geq 0\right) \\
&&\qquad \sim g_{\alpha ,\beta }(0)\mathbf{E}^{+}(z,q)b_{n-j}%
\int_{0}^{K-x}V(-u)du
\end{eqnarray*}%
as $n-j\rightarrow \infty $.

Hence, by (\ref{EstimS_tay}) and the dominated convergence theorem we deduce
that
\begin{eqnarray*}
&&\frac{\mathbf{E}\left[ z^{Z_{n}};\mathcal{R}(n,K),\tau _{n}=j,S_{j}\geq
-N,Z_{j}=q\right] }{b_{n-j}} \\
&&=\int_{-N}^{K\wedge 0}\mathbf{P}\left( S_{j}\in dx,\tau
_{j}=j,Z_{j}=q\right) \\
&&\times \frac{\mathbf{E}\left[ F_{0,n-j}^{q}(z)-F_{0,n-j}^{q}(0);S_{n-j}%
\leq K-x,L_{n-j}\geq 0\right] }{b_{n-j}} \\
& &\sim g_{\alpha ,\beta }(0)\mathbf{E}^{+}(z,q)\int_{-N}^{K\wedge 0}\mathbf{%
P}\left( S_{j}\in dx,\tau _{j}=j,Z_{j}=q\right) \int_{0}^{K-x}V(-u)du
\end{eqnarray*}%
as $n\rightarrow \infty $.

Combining this relation with estimates (\ref{IntermSmall}), (\ref{left1})
and (\ref{left2}), letting $Q$ and $N$ tend to infinity and using (\ref{Defb}%
), we obtain that, for any fixed $j\in \lbrack 0,J]$%
\begin{eqnarray}
&&m_{j}(z):=\lim_{n\rightarrow \infty }\frac{\mathbf{E}\left[ z^{Z_{n}};%
\mathcal{R}(n,K),\tau _{n}=j\right] }{b_{n}}  \notag \\
&&\quad=g_{\alpha ,\beta }(0)\mathbf{E}\left[ \mathbf{E}^{+}(z,Z_{j})%
\int_{0}^{K-S_{j}}V(-u)du\mathbf{;}S_{j}\leq K\wedge 0,\tau _{j}=j\right].
\label{m_j1}
\end{eqnarray}%
Now, letting sequentially $n$ and $J$ tend to infinity we see that%
\begin{equation}
\hat{G}_{left}(K,z):=\lim_{J\rightarrow \infty }\lim_{n\rightarrow \infty }%
\frac{\mathbf{E}\left[ z^{Z_{n}};\mathcal{R}(n,K),\tau _{n}\in \lbrack 0,J]%
\right] }{b_{n}}=\sum_{j=0}^{\infty }m_{j}\left( z\right) .  \label{Gleft}
\end{equation}

Observe now that $\hat{G}_{left}(K,z)>0$ for all $z\in (0,1]$. Indeed, if $%
\mathbf{P}\left( Z_{1}=l\right) >0$ for some $l\in \mathbb{N}$, then
\begin{eqnarray*}
&&\hat{G}_{left}(K,z)\geq m_{1}(z) \\
&&\qquad\quad\geq g_{\alpha ,\beta }(0)\mathbf{E}^{+}(z,l)\int_{-\infty
}^{K\wedge 0}\mathbf{P}\left( Z_{1}=l,S_{1}\in dx\right)
\int_{0}^{K-x}V(-u)du>0.
\end{eqnarray*}

To show that $\hat{G}_{left}(K,z)<\infty $ we introduce the event
\begin{equation*}
A_{u.s}^{(q)}:=\left\{ Z_{0}=q,Z_{n}>0\text{ for all }n\geq 1\right\}
\end{equation*}%
and write an explicit expression for the constant $G_{left}(K)$ in (\ref%
{Limit_surviv}) (see \cite{VDD2024} with a slightly different notation):
\begin{equation}
G_{left}(K)=g_{\alpha ,\beta }(0)\sum_{j=0}^{\infty }\mathbf{E}\left[
\mathbf{P}^{+}\left( A_{u.s}^{(Z_{j})}\right)
\int_{0}^{K-S_{j}}V(-u)du;S_{j}\leq K\wedge 0,\tau _{j}=j\right] .
\label{Gleft0}
\end{equation}%
It was shown in \cite{VDD2024} that $G_{left}(K)<\infty $ which, in turn,
gives $\hat{G}_{left}(K,z)\leq G_{left}(K)<\infty $.

Moreover, using the relation
\begin{equation*}
\lim_{z\uparrow 1}\mathbf{E}^{+}(z,q)=\mathbf{P}^{+}\left(
A_{u.s}^{(q)}\right) ,\ q\geq 1,
\end{equation*}%
and applying two times the dominated convergence theorem to the integrals
and sums we get by (\ref{m_j1}) and (\ref{Gleft}) that%
\begin{equation}
\lim_{z\uparrow 1}\hat{G}_{left}(K,z)=G_{left}(K).  \label{ConvG_left}
\end{equation}

We now assume that $n-j=v\in \lbrack 0,J]\ $and investigate the quantity
\begin{eqnarray*}
\Lambda \left( n-J+1,n\right)  &=&\mathbf{E}[F_{0,n}\left( z\right)
-F_{0,n}\left( 0\right) ;S_{n}\leq K,\tau _{n}\in \lbrack n-J+1,n]] \\
&=&\mathbf{E}[1-F_{0,n}\left( 0\right) ;S_{n}\leq K,\tau _{n}\in \lbrack
n-J+1,n]] \\
&&\,-\mathbf{E}[1-F_{0,n}\left( z\right) ;S_{n}\leq K,\tau _{n}\in \lbrack
n-J+1,n]].
\end{eqnarray*}%
Consider the expectation
\begin{equation*}
\mathbf{E}[1-F_{0,n}\left( z\right) ;S_{n}\leq K,\tau _{n}\in \lbrack
n-J+1,n]].
\end{equation*}%
Using the independency of $F_{1},\ldots ,F_{j}$ and $F_{j+1},\ldots ,F_{n}$,
we write
\begin{eqnarray*}
&&\mathbf{E}[1-F_{0,n}\left( z\right) ;S_{n}\leq K,\tau _{n}=j] \\
&=&\mathbf{E}\left[ 1-F_{0,j}(F_{j,n}(z));S_{j}\leq K-(S_{n}-S_{j}),\tau
_{n}=j\right]  \\
&=&\int_{0}^{\infty }\int_{0}^{1}\mathbf{P}\left( \hat{F}_{0,v}(z)\in du,%
\hat{S}_{v}\in dx,\hat{L}_{v}\geq 0\right) \mathbf{E}\left[
1-F_{0,j}(u);S_{j}\leq K-x,\tau _{j}=j\right] .
\end{eqnarray*}%
To avoid confusion, we temporary use hats in this formula to mark the
respective characteristics of a new branching process generated by the
offspring generating functions $\hat{F}_{1},...,\hat{F}_{v}$ being
independent probabilistic copies of the sequence $F_{1},...,F_{v}$ and
obtained for $n=j+v$ by the following equality in distribution:%
\begin{equation*}
\left( F_{j+1},...,F_{j+v},S_{j+1}-S_{j},...,S_{n}-S_{j}\right) \overset{d}{=%
}\left( \hat{F}_{1},...,\hat{F}_{v},\hat{S}_{1},...,\hat{S}_{v}\right) .
\end{equation*}

It is shown by formulas (7.10)-- (7.12) in \cite{VDD2024} that (see
definition (\ref{Def_O}))
\begin{eqnarray}
\lim_{j\rightarrow \infty }O_{j}\left( u,K-x\right)  &=&\lim_{j\rightarrow
\infty }\frac{\mathbf{E}\left[ 1-F_{0,j}(u);S_{j}\leq K-x,\tau _{j}=j\right]
}{\mathbf{E}\left[ e^{S_{j}};\tau _{j}=j\right] }  \notag \\
&=&h(u,K-x)  \label{Posit222aaa}
\end{eqnarray}%
exists (see formula (\ref{Defh}) below for the precise expression for $%
h(u,K-x)$), where \ $h(u,K-x)>0$ and
\begin{equation*}
\sup_{u\in \lbrack 0,1],w\in (-\infty ,+\infty )}h(u,w)\leq 1
\end{equation*}%
in view of evident inequality $1-F_{0,j}(u)\leq e^{S_{j}}$.

It is easy to deduce for (\ref{Posit222aaa}) and Corollary \ref{C_h} that
there is a constant $C_{1}>0$ such that%
\begin{equation}
h(u,K-x)\leq C_{1}e^{K/2}\left( 1-u\right) e^{-x/2}  \label{Bound_h}
\end{equation}%
for all $0\leq u\leq 1$. Since%
\begin{eqnarray*}
&&\int_{0}^{\infty }\int_{0}^{1}\mathbf{P}\left( F_{0,v}(z)\in du,S_{v}\in
dx,L_{v}\geq 0\right) \left( 1-u\right) e^{-x/2} \\
&&\qquad \qquad \leq \int_{0}^{\infty }\mathbf{P}\left( S_{v}\in
dx,L_{v}\geq 0\right) e^{-x/2}<\infty ,
\end{eqnarray*}%
we conclude by (\ref{Posit222aaa}) and the dominated convergence theorem
that
\begin{eqnarray}
&&\lim_{j\rightarrow \infty }\int_{0}^{\infty }\int_{0}^{1}\mathbf{P}\left(
F_{0,v}(z)\in du,S_{v}\in dx,L_{v}\geq 0\right) O_{j}\left( u,K-x\right)
\notag \\
&&\quad =\int_{0}^{\infty }\int_{0}^{1}\mathbf{P}\left( F_{0,v}(z)\in
du,S_{v}\in dx,L_{v}\geq 0\right) h(u,K-x)  \notag \\
&&\qquad \qquad =\mathbf{E}\left[ h(F_{0,v}(z),K-S_{v})\mathbf{I}\left\{
L_{v}\geq 0\right\} \right] .  \label{SingleTermRight}
\end{eqnarray}%
Further we have%
\begin{eqnarray*}
&&\sum_{v=0}^{\infty }\int_{0}^{\infty }\int_{0}^{1}\mathbf{P}\left(
F_{0,v}(z)\in du,S_{v}\in dx,L_{v}\geq 0\right) O_{j}\left( u,K-x\right)  \\
&&\quad \leq C\sum_{v=0}^{\infty }\int_{0}^{\infty }\int_{0}^{1}\mathbf{P}%
\left( F_{0,v}(z)\in du,S_{v}\in dx,L_{v}\geq 0\right) \left( 1-u\right)
e^{-x/2} \\
&&\quad \leq C\sum_{v=0}^{\infty }\int_{0}^{\infty }\mathbf{P}\left(
S_{v}\in dx,L_{v}\geq 0\right) e^{-x/2}=C\int_{0}^{\infty
}e^{-x/2}U(dx)<\infty .
\end{eqnarray*}%
Hence, applying two times the dominated convergence theorem and taking into
account that if $j=n-v\sim n$ as $n\rightarrow \infty $ then
\begin{equation*}
\mathbf{E}\left[ e^{S_{j}};\tau _{j}=j\right] \sim \mathbf{E}\left[
e^{S_{n}};\tau _{n}=n\right]
\end{equation*}%
according to (\ref{AsymExponent}), we get%
\begin{eqnarray*}
&&\lim_{J\rightarrow \infty }\lim_{n\rightarrow \infty }\frac{\mathbf{E}%
[1-F_{0,n}\left( z\right) ;S_{n}\leq K,\tau _{n}\in \lbrack n-J+1,n]]}{%
\mathbf{E}\left[ e^{S_{n}};\tau _{n}=n\right] } \\
&=&\lim_{J\rightarrow \infty }\lim_{n\rightarrow \infty }\sum_{v=0}^{J}\frac{%
\mathbf{E}[1-F_{0,n}\left( z\right) ;S_{n}\leq K,\tau _{n}=n-v]]}{\mathbf{E}%
\left[ e^{S_{n-v}};\tau _{n-v}=n-v\right] }\frac{\mathbf{E}\left[
e^{S_{n-v}};\tau _{n-v}=n-v\right] }{\mathbf{E}\left[ e^{S_{n}};\tau _{n}=n%
\right] } \\
&=&\lim_{J\rightarrow \infty }\lim_{n\rightarrow \infty
}\sum_{v=0}^{J}\int_{0}^{\infty }\int_{0}^{1}\mathbf{P}\left( F_{0,v}(z)\in
du,S_{v}\in dx,L_{v}\geq 0\right) O_{n-v}\left( u,K-x\right)  \\
&&\qquad =\lim_{J\rightarrow \infty }\sum_{v=0}^{J}\mathbf{E}\left[ \mathbf{I%
}\left\{ L_{v}\geq 0\right\} \lim_{j\rightarrow \infty }O_{j}\left(
F_{0,v}(z),K-S_{v}\right) \right]  \\
&&\qquad =\mathbf{E}\left[ \sum_{v=0}^{\infty }h(F_{0,v}(z),K-S_{v})\mathbf{I%
}\left\{ L_{v}\geq 0\right\} \right] <\infty .
\end{eqnarray*}%
Taking this fact into account, using (\ref{AsymS_tay}) and summing (\ref%
{SingleTermRight}) over $v$ from $0$ to $\infty $ we get
\begin{eqnarray}
&&\lim_{J\rightarrow \infty }\lim_{n\rightarrow \infty }\frac{\mathbf{E}%
[F_{0,n}\left( z\right) -F_{0,n}\left( 0\right) ;S_{n}\leq K,\tau _{n}\in
\lbrack n-J+1,n]]}{b_{n}}  \notag \\
&&\qquad \qquad \qquad =\hat{G}_{right}(K,z):=G_{right}(K)-G_{right}(K,z),
\label{NEWGRIGH0}
\end{eqnarray}%
where
\begin{eqnarray}
&&G_{right}(K)=g_{\alpha ,\beta }(0)\int_{0}^{\infty }e^{-y}U(y)dy  \notag \\
&&\times \mathbf{E}\left[ \sum_{v=0}^{\infty }h(F_{0,v}(0),K-S_{v})\mathbf{I}%
\left\{ L_{v}\geq 0\right\} \right]   \label{Gright0}
\end{eqnarray}%
is the same constant as in (\ref{Limit_surviv})\ (see formula (7.14) in \cite%
{VDD2024} with a slightly different notation) and%
\begin{eqnarray*}
&&G_{right}(K,z):=g_{\alpha ,\beta }(0)\int_{0}^{\infty }e^{-y}U(y)dy \\
&&\times \mathbf{E}\left[ \sum_{v=0}^{\infty }h(F_{0,v}(z),K-S_{v})\mathbf{I}%
\left\{ L_{v}\geq 0\right\} \right] .
\end{eqnarray*}%
Combining (\ref{IntermSmall}) with (\ref{Gleft}) and (\ref{NEWGRIGH0}) we
obtain
\begin{equation*}
\lim_{n\rightarrow \infty }\frac{\mathbf{E}\left[ z^{Z_{n}};\mathcal{R}(n,K)%
\right] }{b_{n}}=\hat{G}_{left}(K,z)+\hat{G}_{right}(K,z).
\end{equation*}%
Taking (\ref{Limit_surviv}) into account we conclude that%
\begin{equation}
\lim_{n\rightarrow \infty }\mathbf{E}\left[ z^{Z_{n}}|\mathcal{R}(n,K)\right]
=\frac{\hat{G}_{left}(K,z)+\hat{G}_{right}(K,z)}{G_{left}(K)+G_{right}(K)}.
\label{DiscreteLimit}
\end{equation}

Now we show that the limiting discrete distribution (\ref{DiscreteLimit}) is
proper. Since $\lim_{z\uparrow 1}\hat{G}_{left}(K,z)=G_{left}(K)$ (see (\ref%
{ConvG_left})), it remains to check that%
\begin{equation*}
\lim_{z\uparrow 1}\hat{G}_{right}(K,z)=G_{right}(K).
\end{equation*}%
To this aim we use (\ref{Bound_h}) and the estimates%
\begin{eqnarray*}
G_{right}(K,z) &\leq &C\mathbf{E}\left[ \sum_{v=0}^{\infty }\left(
1-F_{0,v}(z)\right) e^{-S_{v}/2}\mathbf{I}\left\{ L_{v}\geq 0\right\} \right]
\\
&\leq &R_{1}(\varepsilon ,z)+R_{2}(\varepsilon ,z),
\end{eqnarray*}%
where%
\begin{equation*}
R_{1}(\varepsilon ,z)=C\varepsilon \mathbf{E}\left[ \sum_{v=0}^{\infty
}e^{-S_{v}/2}\mathbf{I}\left\{ 1-F_{0,v}(z)\leq \varepsilon ,L_{v}\geq
0\right\} \right]
\end{equation*}%
and%
\begin{equation*}
R_{2}(\varepsilon ,z)=C\mathbf{E}\left[ \sum_{v=0}^{\infty }e^{-S_{v}/2}%
\mathbf{I}\left\{ 1-F_{0,v}(z)\geq \varepsilon ,L_{v}\geq 0\right\} \right] .
\end{equation*}%
Clearly,%
\begin{eqnarray}
R_{1}(\varepsilon ,z) &\leq &C\varepsilon \mathbf{E}\left[
\sum_{v=0}^{\infty }e^{-S_{v}/2}\mathbf{I}\left\{ L_{v}\geq 0\right\} \right]
\notag \\
&=&C\varepsilon \int_{0}^{\infty }e^{-x/2}U(dx)\leq C_{1}\varepsilon .
\label{EstR_1}
\end{eqnarray}%
To evaluate $R_{2}(\varepsilon ,z)$ observe that, for any $\varepsilon _{1}>0
$ there exists $M=M(\varepsilon _{1})$ such that
\begin{eqnarray}
&&\mathbf{E}\left[ \sum_{v=0}^{\infty }e^{-S_{v}/2}\mathbf{I}\left\{
1-F_{0,v}(z)\geq \varepsilon ,S_{v}\geq M,L_{v}\geq 0\right\} \right]
\notag \\
&&\qquad \qquad \qquad \leq \mathbf{E}\left[ \sum_{v=0}^{\infty }e^{-S_{v}/2}%
\mathbf{I}\left\{ S_{v}\geq M,L_{v}\geq 0\right\} \right]   \notag \\
&&\qquad \qquad \qquad \qquad =\int_{M}^{\infty }e^{-x/2}U(dx)<\varepsilon
_{1}.  \label{EstR_23}
\end{eqnarray}%
Further, using (\ref{EstimS_tay}) and observing that%
\begin{equation*}
\sum_{v=0}^{\infty }b_{v}<\infty
\end{equation*}%
in view of (\ref{Defb}) we conclude that, for any $\varepsilon _{1}>0$ there
exists $J=J(\varepsilon _{1},M)$ such that
\begin{eqnarray}
&&\mathbf{E}\left[ \sum_{v=J+1}^{\infty }e^{-S_{v}/2}\mathbf{I}\left\{
1-F_{0,v}(z)\geq \varepsilon ,S_{v}\leq M,L_{v}\geq 0\right\} \right]
\notag \\
&&\qquad \qquad \qquad \qquad \leq \sum_{v=J+1}^{\infty }\mathbf{P}\left(
S_{v}\leq M,L_{v}\geq 0\right)   \notag \\
&&\qquad \qquad \qquad \qquad \leq C\int_{0}^{M}V(-u)du\sum_{v=J+1}^{\infty
}b_{j}<\varepsilon _{1}.  \label{EstR_22}
\end{eqnarray}%
Observe that estimates (\ref{EstR_1})-(\ref{EstR_22}) are uniform in $z\in
\lbrack 0,1]$.

Finally, for any fixed $v$
\begin{equation*}
\lim_{z\uparrow 1}F_{0,v}(z)=F_{0,v}(1)=\mathbf{P}\left( Z_{v}<\infty |%
\mathcal{E}\right) \mathbf{=}1\quad \text{ }\mathbf{P}-\text{a.s.}
\end{equation*}%
Therefore, for any fixed $J$%
\begin{eqnarray*}
&&\lim_{z\uparrow 1}\mathbf{E}\left[ \sum_{v=0}^{J}e^{-S_{v}/2}\mathbf{I}%
\left\{ 1-F_{0,v}(z)\geq \varepsilon ,S_{v}\leq M,L_{v}\geq 0\right\} \right]
\\
&&\qquad \qquad \leq \lim_{z\uparrow 1}\sum_{v=0}^{J}\mathbf{P}\left(
1-F_{0,v}(z)\geq \varepsilon \right) \leq \frac{1}{\varepsilon }%
\sum_{v=0}^{J}\lim_{z\uparrow 1}\mathbf{E}\left[ 1-F_{0,v}(z)\right] \\
&&\qquad \qquad\qquad =\frac{1}{\varepsilon }\sum_{v=0}^{J}\mathbf{E}\left[
1-\lim_{z\uparrow 1}F_{0,v}(z)\right] =0.
\end{eqnarray*}

Hence it follows that
\begin{equation*}
\lim_{z\uparrow 1}\hat{G}_{right}(K,z)=G_{right}(K)-\lim_{z\uparrow
1}G_{right}(K,z)=G_{right}(K).
\end{equation*}%
Combining the estimates above we conclude that%
\begin{equation*}
\lim_{z\uparrow 1}\frac{\hat{G}_{left}(K,z)+\hat{G}_{right}(K,z)}{%
G_{left}(K)+G_{right}(K)}=1
\end{equation*}%
and, therefore, the conditional laws $\mathcal{L}(Z_{n}|Z_{n}>0;S_{n}\leq
K),n\geq 1,$ converge weakly, as $n\rightarrow \infty $ to a proper
probability distribution with support contained in the set of positive
integers.

Theorem \ref{T_simplelaw} is proved.

\section{P\textbf{roof of Theorem \protect\ref{T_generalcase}\label{sec4}}}

In this section we will use from time to time the notation $Z_{1}^{\left(
q\right) },Z_{2}^{\left( q\right) },...$ if $Z_{0}=q\in \mathbb{N}$. We
introduce, for integers $m\geq k\geq j$ the process
\begin{equation}
\mathcal{Y}_{(q)}^{k,m,j}:=\{Y_{(q)}^{k,m,j}(t),\quad t\in \lbrack 0,1]\},
\label{DEfProcY}
\end{equation}%
where\bigskip
\begin{equation}
Y_{(q)}^{k,m,j}(t):=\exp \left\{ -S_{k+\left[ \left( m-2k\right) t\right]
-j}\right\} Z_{k+\lfloor (m-2k)t\rfloor -j}^{(q)}.  \label{DefComponent_Y}
\end{equation}%
In \ words: $\mathcal{Y}_{(q)}^{k,m,j}$ is a piece of the trajectory of the
rescaled population size process
\begin{equation*}
\mathcal{Y}_{(q)}^{0,m}=\left\{ \exp \left\{ -S_{r}\right\}
Z_{r}^{(q)},r=0,1,...,m\right\}
\end{equation*}%
between the moments $k-j$ and $m-k-j$ given that the original process starts
by $q$ individuals at time $0$. \ Note that (see \ref{Def_Y})%
\begin{equation*}
\mathcal{Y}_{(1)}^{k,m,0}=\mathcal{Y}^{k,m}=\{Y^{k,m}(t),\quad t\in \lbrack
0,1]\}
\end{equation*}%
if $j=0$ and $Z_{0}=q=1$.

We need to analyze the process $\mathcal{Y}^{\theta n,n}=\{Y^{\theta
n,n}(t),t\in \lbrack 0,1]\}$, where%
\begin{equation*}
Y^{\theta n,n}(t):=\exp \left\{ -S_{\theta n+\lfloor (1-2\theta )nt\rfloor
}\right\} Z_{\theta n+\lfloor (1-2\theta )nt\rfloor }.
\end{equation*}

Let $\phi $ be a bounded continuous function on $D\left[ 0,1\right] $. We
assume, without loss of generality, that $0\leq \phi \leq 1$ and $\phi
\not\equiv 0.$ Our aim is to show that
\begin{equation*}
\lim_{n\rightarrow \infty }\mathbf{E[}\phi \left( \mathcal{Y}^{\theta
n,n}\right) |\mathcal{R}(n,K)]
\end{equation*}%
exists. We fix $J\in \lbrack 0,n]$ and write the representation
\begin{equation*}
\mathbf{E[}\phi \left( \mathcal{Y}^{\theta n,n}\right) ;\mathcal{R}%
(n,K)]=\Gamma \left( 0,J\right) +\Gamma \left( J+1,n-J\right) +\Gamma \left(
n-J+1,n\right) ,
\end{equation*}%
where for $0\leq n_{1}<n_{2}$%
\begin{eqnarray*}
\Gamma \left( n_{1},n_{2}\right)  &=&\mathbf{E[}\phi \left( \mathcal{Y}%
^{\theta n,n}\right) ;\mathcal{R}(n,K),\tau _{n}\in \lbrack n_{1},n_{2}]] \\
&=&\sum\limits_{j=n_{1}}^{n_{2}}\mathbf{E[}\phi \left( \mathcal{Y}^{\theta
n,n}\right) ;\mathcal{R}(n,K),\tau _{j}=j,\min_{j<l\leq n}S_{l}\geq S_{j}].
\end{eqnarray*}%
First we observe that%
\begin{equation*}
\lim_{J\rightarrow \infty }\lim_{n\rightarrow \infty }\frac{\Gamma \left(
J+1,n-J\right) }{b_{n}}=0
\end{equation*}%
in view of (\ref{441}) and the estimate%
\begin{equation*}
\left\vert \Gamma \left( J+1,n-J\right) \right\vert \leq \mathbf{P}\left(
\mathcal{R}(n,K),\tau _{n}\in \left[ J+1,n-J\right] \right) .
\end{equation*}

Now we study the behavior of $\Gamma \left( 0,J\right) $ for large but fixed
$J$.

Using the definitions (\ref{DEfProcY})-(\ref{DefComponent_Y}) we conclude by
the estimates similar to (\ref{left1}) and (\ref{left2}) that, for any $j\in
\lbrack 0,J]$ and $\theta n\geq J$
\begin{eqnarray*}
&&\mathbf{E[}\phi \left( \mathcal{Y}^{\theta n,n}\right) ;\mathcal{R}%
(n,K),\tau _{n}=j]=\varepsilon _{j,N,Q}b_{n-j} \\
&&+\mathbf{E[}\phi \left( \mathcal{Y}^{\theta n,n}\right) ;\mathcal{R}%
(n,K),\tau _{n}=j,S_{j}>-N,Z_{j}\leq Q] \\
&=&\varepsilon _{j,N,Q}b_{n-j} \\
&&+\int_{-N}^{0}\sum_{q=1}^{Q}\mathbf{P}\left( S_{j}\in dx,Z_{j}=q,\tau
_{j}=j\right) \\
&&\times \mathbf{E}\left[ \phi \left( e^{-x}\mathcal{\hat{Y}}_{(q)}^{n\theta
,n,j}\right) \mathbf{I}\left\{ \hat{Z}_{n-j}^{(q)}>0\right\} ;\hat{S}%
_{n-j}\leq K-x,\hat{L}_{n-j}\geq 0\right] ,
\end{eqnarray*}%
where%
\begin{equation*}
\limsup_{N\rightarrow \infty }\limsup_{Q\rightarrow \infty }\sup_{0\leq
j\leq J}\left\vert \varepsilon _{j,N,Q}\right\vert =0
\end{equation*}%
and hats are temporary used to denote the respective characteristics of a
new branching process $\left\{ \hat{Z}_{m},m=0,1,...\right\} $ in random
environment which is independent of $\left\{ Z_{r},r=0,1,...,j\right\} $ and
has the same distribution as the initial process.

For $x\in \left( -\infty ,+\infty \right) $ and $q\in \mathbb{N}$ let $%
\mathcal{W}_{\left( q\right) }^{x}=\left\{ W_{\left( q\right) }^{x}(t),t\in
\lbrack 0,1]\right\} $ be the process with constant paths
\begin{equation*}
W_{\left( q\right) }^{x}(t):=e^{-x}W_{\left( q\right) }^{+},t\in \lbrack
0,1],
\end{equation*}%
where the random variable $W_{\left( q\right) }^{+}$ is specified by (\ref%
{DefWPlus}). Set, for brevity,
\begin{equation*}
\mathcal{W}_{\left( q\right) }:=\mathcal{W}_{\left( q\right) }^{0},\quad
\mathcal{W}^{x}:=\mathcal{W}_{\left( 1\right) }^{x}
\end{equation*}%
and
\begin{equation*}
\quad \mathcal{W}=\mathcal{W}^{0}:=\left\{ W(t),t\in \lbrack 0,1]\right\}
=\left\{ W^{+},t\in \lbrack 0,1]\right\} .
\end{equation*}%
Let
\begin{equation*}
H_{q,\infty }\left( x\right) :=\phi \left( \mathcal{W}_{\left( q\right)
}^{x}\right) \mathbf{I}\left\{ W_{\left( q\right) }^{+}>0\right\}
\end{equation*}%
and
\begin{equation*}
D_{j}\left( \phi ,Q\right) :=\mathbf{E}\left[ \mathbf{E}^{+}\left[
H_{Z_{j},\infty }\left( S_{j}\right) \right] \mathbf{I}\left\{ \tau
_{j}=j,Z_{j}\in \lbrack 1,Q]\right\} \right] ,\quad Q\in \mathbb{N}.
\end{equation*}

Since $0\leq \phi \leq 1$, it follows that
\begin{equation*}
\mathbf{E}\left[ \mathbf{E}^{+}\left[ H_{q,\infty }\left( S_{j}\right) %
\right] \mathbf{I}\left\{ \tau _{j}=j,Z_{j}=q\right\} \right] \leq \mathbf{P}%
\left( Z_{j}=q,\tau _{j}=j\right)
\end{equation*}%
and%
\begin{equation*}
\sum_{q=1}^{\infty }\mathbf{P}\left( Z_{j}=q,\tau _{j}=j\right) \leq \mathbf{%
P}\left( \tau _{j}=j\right) \leq 1.
\end{equation*}%
Hence, by monotonicity we conclude that
\begin{equation*}
\lim_{Q\rightarrow \infty }D_{j}\left( \phi ,Q\right) =:D_{j}\left( \phi
\right) <\infty
\end{equation*}%
exists, where%
\begin{equation*}
D_{j}\left( \phi \right) :=\mathbf{E}\left[ \mathbf{E}^{+}\left[
H_{Z_{j},\infty }\left( S_{j}\right) \right] \mathbf{I}\left\{ \tau
_{j}=j,Z_{j}>0\right\} \right] .
\end{equation*}

According to Lemma \ref{Lemprop3.1}, for any $q\in \mathbb{N}$ the process $%
e^{-x}\mathcal{\hat{Y}}_{\left( q\right) }^{\theta n,n,j}$ converges $%
\mathbf{P}^{+}$-a.s. as $n\rightarrow \infty $ to the process $\mathcal{W}%
_{\left( q\right) }^{x}$ in the uniform metric and, therefore, in the
Skorokhod metric in $D[0,1]$. Thus, for any $q\in \mathbb{N}$
\begin{eqnarray*}
H_{q,\theta ,n,j}\left( x\right) := &&\phi \left( e^{-x}\mathcal{\hat{Y}}%
_{\left( q\right) }^{\theta n,n,j}\right) \mathbf{I}\left\{ \hat{Z}%
_{n-j}^{\left( q\right) }>0\right\}  \\
&\rightarrow &\phi \left( \mathcal{W}_{\left( q\right) }^{x}\right) \mathbf{I%
}\left\{ W_{\left( q\right) }^{+}>0\right\} =:H_{q,\infty }\left( x\right)
\end{eqnarray*}%
$\mathbf{P}^{+}$-a.s. as $n\rightarrow \infty $. Using now Lemma \ref{L_cond}
we see that, as $n\rightarrow \infty $
\begin{eqnarray*}
\mathbf{E}\left[ H_{q,\theta ,n,j}\left( x\right) ;S_{n-j}\leq
K-x,L_{n-j}\geq 0\right]  &=&\left( \mathbf{E}^{+}\left[ H_{q,\infty }\left(
x\right) \right] +o\left( 1\right) \right)  \\
&\times &\mathbf{P}\left( S_{n-j}\leq K-x,L_{n-j}\geq 0\right) .
\end{eqnarray*}%
Assume that%
\begin{equation*}
\mathbf{E}^{+}\left[ H_{q,\infty }\left( x\right) \right] >0
\end{equation*}%
(if this is not the case then the arguments to follow are simplified).
Taking into account (\ref{AsymS_tay}) we conclude that, given Condition $B0$%
, for each $x\in \lbrack -N,K]$ and any $q\in \lbrack 1,Q]$%
\begin{eqnarray*}
&&\mathbf{E}\left[ \phi \left( e^{-x}\mathcal{\hat{Y}}_{\left( q\right)
}^{\theta n,n,j}\right) \mathbf{I}\left\{ \hat{Z}_{n-j}^{\left( q\right)
}>0\right\} ;\hat{S}_{n-j}\leq K-x,\hat{L}_{n-j}\geq 0\right]  \\
&&\qquad \qquad \qquad \qquad \sim \mathbf{E}^{+}\left[ H_{q,\infty }\left(
x\right) \right] g_{\alpha ,\beta }(0)b_{n}\int_{0}^{K-x}V\left( -u\right) du
\end{eqnarray*}%
as $n\rightarrow \infty $.

Using (\ref{EstimS_tay}) we obtain by the dominated convergence theorem that
\begin{eqnarray*}
&&\frac{\mathbf{E[}\phi \left( \mathcal{Y}^{\theta n,n}\right) ;\tau _{n}=j,%
\mathcal{R}(n,K),S_{j}\geq -N,Z_{j}=q]}{b_{n-j}} \\
&=&\int_{-N}^{K\symbol{94}0}\mathbf{P}\left( S_{j}\in dx,Z_{j}=q,\tau
_{j}=j\right) \\
&&\times \frac{\mathbf{E}\left[ \phi \left( e^{-x}\mathcal{\hat{Y}}_{\left(
q\right) }^{\theta n,n,j}\right) \mathbf{I}\left\{ \hat{Z}_{n-j}^{\left(
q\right) }>0\right\} ;\hat{S}_{n-j}\leq K-x,\hat{L}_{n-j}\geq 0\right] }{%
b_{n-j}} \\
&\thicksim &g_{\alpha ,\beta }\left( 0\right) \mathbf{E}\left[ \mathbf{E}^{+}%
\left[ H_{Z_{j},\infty }\left( S_{j}\right) \right] \int_{0}^{K-S_{j}}V%
\left( -u\right) du;-N\leq S_{j}\leq K\wedge 0,\tau _{j}=j\right]
\end{eqnarray*}%
as $n\rightarrow \infty .$

Combining all these estimates, letting $Q$ and $N$ tend to infinity, using~
the equivalence $b_{n}\sim b_{n-j},n\rightarrow \infty ,$ valid for each
fixed $j$ in view of (\ref{Defb}), we deduce that, for any fixed $j\in
\lbrack 0,J]$%
\begin{eqnarray*}
&&m_{j}\left(K,\phi \right):=\lim_{n\rightarrow \infty }\frac{\mathbf{E}%
[\phi \left( \mathcal{Y}^{\theta n,n}\right) ;\mathcal{R}(n,K),\tau _{n}=j]}{%
b_{n}} \\
&&\qquad=g_{\alpha,\beta}(0)\mathbf{E}\left[\mathbf{E}^{+}\left[
H_{Z_{j},\infty}\left( S_{j}\right) \right]\int_{0}^{K-S_{j}}V\left(
-u\right)du; S_{j}\leq K\wedge0,\tau _{j}=j\right].
\end{eqnarray*}

Thus,%
\begin{eqnarray*}
G_{left}\left( K,\phi \right)&:=&\lim_{J\rightarrow \infty
}\lim_{n\rightarrow \infty }\frac{\mathbf{E}[\phi \left( \mathcal{Y}^{\theta
n,n}\right) ;\mathcal{R}(n,K),\tau _{n}\in \lbrack 0,J]]}{b_{n}} \\
&=&\sum_{j=0}^{\infty }m_{j}\left( K,\phi \right) =\int \phi (w)\lambda
_{left}(dw),
\end{eqnarray*}%
where $\lambda _{left}(\cdot )$ is the measure on the space $D[0,1]\ $given
by%
\begin{equation*}
\lambda _{left}(dw):=g_{\alpha ,\beta }(0)\sum_{j=0}^{\infty }\mathbf{E}%
\left[ \lambda _{Z_{j},S_{j}}(dw);S_{j}\leq K\wedge 0,\tau _{j}=j\right]
\end{equation*}%
with%
\begin{equation*}
\lambda _{q,x}(dw):=\mathbf{P}^{+}\left( W_{(q)}^{x}\in
dw,W_{(q)}^{+}>0\right) \int_{0}^{K-x}V\left( -u\right) du.
\end{equation*}%
In particular,%
\begin{equation}
\int \lambda _{left}(dw)=G_{left}(K).  \label{massLeft}
\end{equation}%
By Lemma \ref{Lemprop3.1} the total mass of $\lambda _{q,x}$ is
\begin{equation*}
\mathbf{P}^{+}(Z_{0}=q,Z_{n}>0\text{ for all }n)\int_{0}^{K-x}V\left(
-u\right) du
\end{equation*}%
and $\lambda _{q,x}$ puts its entire mass on strictly positive constant
functions, and, hence, so does $\lambda _{left}.$

Finiteness of $G_{left}\left( K,\phi \right) $ is a corollary of the
inequality%
\begin{equation*}
\frac{\mathbf{E}[\phi \left( \mathcal{Y}^{\theta n,n}\right) ;\mathcal{R}%
(n,K),\tau _{n}\in \lbrack 0,J]]}{b_{n}}\leq \frac{\mathbf{P}\left( \mathcal{%
R}(n,K)\right) }{b_{n}}
\end{equation*}%
and (\ref{Limit_surviv}).

Consider now the sum $\Gamma \left( n-J+1,n\right) $ and write $%
n-j=v\in \lbrack 0,J],$ i.e. $n=j+v$.

Let $v$ be fixed. As above we assume, without loss of generality, that $\phi
$ is a continuous function on $D\left[ 0,1\right] $ meeting the inequalities
$0\leq \phi \left( \cdot \right) \leq 1$ and $\phi \not\equiv 0$.$~$\ Using
the independency of $F_{1},\ldots ,F_{j}$ from $F_{j+1},\ldots ,F_{n}$, we
write for $(1-\theta )n\leq j$
\begin{eqnarray*}
&&\mathbf{E}\left[ \phi \left( \mathcal{Y}^{\theta n,n}\right) ;\mathcal{R}%
(n,K),\tau _{n}=j\right] \\
&&\quad =\mathbf{E}\left[ \mathbf{E}\left[ \phi \left( \mathcal{Y}^{\theta
n,n}\right) ;\mathcal{R}(n,K),\tau _{n}=j|\mathcal{E}\right] \right] \\
&&\quad =\mathbf{E}\left[ \phi \left( \mathcal{Y}^{\theta n,n}\right) \left(
1-F_{(1-\theta )n,n}^{Z_{(1-\theta )n}}(0)\right) ;S_{n}\leq K,\tau _{n}=j%
\right] \\
&&\quad =\mathbf{E}\left[ \phi \left( \mathcal{Y}^{\theta n,n}\right) \left(
1-F_{(1-\theta )n,j}^{Z_{(1-\theta )n}}(F_{j,n}(0))\right) ;S_{j}\leq
K-(S_{n}-S_{j}),\tau _{n}=j\right] \\
&&\quad =\int_{0}^{\infty }\int_{0}^{1}\mathbf{P}\left( \hat{F}_{0,v}(0)\in
du,\hat{S}_{v}\in dx,\hat{L}_{v}\geq 0\right) \\
&&\quad \qquad \times \mathbf{E}\left[ \phi \left( \mathcal{Y}^{\theta
n,n}\right) \left( 1-F_{(1-\theta )n,j}^{Z_{(1-\theta )n}}(u)\right)
;S_{j}\leq K-x,\tau _{j}=j\right] .
\end{eqnarray*}

Our aim is to investigate, as $n\rightarrow \infty $ the asymptotic behavior
of the quantity
\begin{equation*}
\mathbf{E}\left[ \phi \left( \mathcal{Y}^{\theta n,n}\right) \left(
1-F_{(1-\theta )n,j}^{Z_{(1-\theta )n}}(u)\right) ;S_{j}\leq K-x,\tau _{j}=j%
\right] \qquad
\end{equation*}%
for $j\in \lbrack n-J,n]$, fixed $u\in \lbrack 0,1)$ and $x\geq K$.

In view of the estimate
\begin{equation*}
(1-\theta )n=(1-\theta )\left( j+v\right) \leq (1-\theta )\left( j+J\right)
\end{equation*}%
there exists $\delta \in \left( 1-\theta ,1\right) $ such that $\left(
1-\theta \right) \left( j+J\right) <\delta j$ for sufficiently large $j\in
\lbrack n-J,n].$

According to Corollary \ref{C_polish} the space $\mathcal{A}:=D[0,1]\times
\mathbb{R}^{+}$ is polish with norm $d_{0}\times \left\Vert x\right\Vert
_{\infty }$. Keeping in mind this fact, using relation $n=j+v$ and assuming
that $v\in \lbrack 0,J]$ is fixed we consider the sequence of elements
\begin{eqnarray*}
\mathbf{W}_{j}:= &&\left( W_{j0},W_{j1}\right) =\left( \mathcal{Y}^{\theta
(j+v),j+v},Y^{\theta (j+v),j+v}(1)\right)  \\
&=&\left( \mathcal{Y}^{\theta n,n},Y^{\theta n,n}(1)\right) ,\;n=j+v\geq 1,
\end{eqnarray*}%
of the polish space. Clearly, this sequence may be considered as a sequence
of the form
\begin{equation*}
\mathfrak{a}_{j}\left( F_{0},...,F_{\lfloor \delta j\rfloor
};Z_{0},...,Z_{\lfloor \delta j\rfloor }\right) ,\;j\geq 1,
\end{equation*}%
taking values in $\mathcal{A}$.

Introduce the functions%
\begin{equation*}
\Psi _{K-x}(w_{1},b,y):=\left( 1-b^{w_{1}\exp \left\{ y\right\} }\right)
e^{-y}\mathbf{I}\left\{ y\leq K-x\right\}
\end{equation*}%
and
\begin{equation*}
\vec{\Psi}_{K-x}(\mathbf{w},b,y):=\phi \left( w\right) \Psi _{K-x}(w_{1},b,y)
\end{equation*}%
for $\mathbf{w}=\left( w,w_{1}\right) \in $ $\mathcal{A}$, $0\leq b\leq
1,y\leq 0$ with $0^{0}=1,$ and continue $\Psi _{K-x}$ to the other values of
$w_{1},b,y\leq 0$ to a bounded smooth function. In doing so, discontinuities
in points where $b=0$ and in points where $y=K-x$ are unavoidable.

Setting
\begin{equation*}
\tilde{B}_{j}(u):=\left( F_{(1-\theta )n,j}(u)\right) ^{\exp \left\{
S_{(1-\theta )n}-S_{j}\right\} }
\end{equation*}%
and using the representation%
\begin{equation*}
1-F_{(1-\theta )n,j}^{Z_{(1-\theta )n}}(u)=\left( 1-\left( \tilde{B}%
_{j}(u)\right) ^{W_{j1}\exp \left\{ S_{j}\right\} }\right)
e^{-S_{j}}e^{S_{j}}
\end{equation*}%
we write the identity
\begin{eqnarray}
&&\mathbf{E}\left[ \phi \left( \mathcal{Y}^{\theta n,n}\right) \left(
1-F_{(1-\theta )n,j}^{Z_{(1-\theta )n}}(u)\right) \mathbf{I}\left\{
S_{j}\leq K-x\right\} ;\tau _{j}=j\right]  \notag \\
&&\qquad \qquad \qquad \qquad \qquad =\mathbf{E}\left[ \vec{\Psi}%
_{K-x}\left( \mathbf{W}_{j};\tilde{B}_{j}(u),S_{j}\right) e^{S_{j}};\tau
_{j}=j\right] .  \label{Identit}
\end{eqnarray}%
Our aim is to apply Lemma \ref{L_Subcr} to the function $\vec{\Psi}_{K-x}$.
However, one use this lemma only to the bounded continuous functions. We
show that in the case of BPRE's this difficulty can be bypassed.

It follows from Lemma \ref{Lemprop3.1} that if $n\rightarrow \infty $ and $%
\theta \in (0,1/2),$ then the process $\mathcal{Y}^{\theta n,n}=\mathcal{Y}%
^{\theta (j+v),(j+v)(1-\theta )}$ converges, as $j\rightarrow \infty $ to
the process
\begin{equation*}
\mathcal{W}=\left\{ W(t),\ 0\leq t\leq 1\right\} =\left\{ W^{+},\ 0\leq
t\leq 1\right\}
\end{equation*}%
$\mathbf{P}_{y}^{+}$- a.s. for all $y\geq 0$ in the uniform metric and,
therefore, in the Skorokhod metric on $D[0,1]$. Besides, the limit is one
and the same for all $v\in \lbrack 0,J]$. Using this fact we conclude that,
as $j\rightarrow \infty $
\begin{equation*}
\mathbf{W}_{j}=\left( W_{j,0},W_{j,1}\right) \rightarrow \left( \mathcal{W}%
,W^{+}\right) :=\mathbf{W}_{\infty }
\end{equation*}%
$\mathbf{P}_{y}^{+}$- a.s. for all $y\geq 0$.

Further, it follows from Lemma 3.2 in~\cite{ABGV2011}\ that the sequence
\begin{eqnarray*}
B_{j}(u)&:=&\left( F_{j-(1-\theta )n,0}(u)\right) ^{\exp \left\{
-S_{j-(1-\theta )n}\right\} } \\
&=&\left( F_{n\theta -v,0}(u)\right) ^{\exp \left\{ -S_{n\theta -v}\right\}
},\,j=1,2,...
\end{eqnarray*}%
is nondecreasing in $j~$\ and $B_{j}(u)\geq u^{\exp \left\{ -S_{0}\right\} }$%
. Thus, as $j\rightarrow \infty $
\begin{equation*}
B_{j}(u)\rightarrow B_{\infty }(u)\geq u^{\exp \left\{ -S_{0}\right\}
}=u\qquad \mathbf{P}^{-}\text{- a.s.}
\end{equation*}%
In view of
\begin{equation*}
\mathbf{P}\left( \tilde{B}_{j}(u)\geq u\right) =\mathbf{P}\left(
B_{j}(u)\geq u\right) =1,
\end{equation*}%
it follows that the second random term in $\vec{\Psi}_{K-x}(\mathbf{W}_{j},%
\tilde{B}_{j}(u),S_{j})$ is separated from zero with probability one.

Therefore, aiming to apply Lemma \ref{L_Subcr} to the function $\vec{\Psi}%
_{K-x}(\mathbf{w},b,y)$ we can consider this function in the domain,\ where $%
y<K-x,$ only\ for $b\in (0,1].$ Clearly, $\Psi _{K-x}$ is continuous in this
domain. Finally, the discontinuity in the point $y=K-x$ has probability $0$
with respect to the measure $\mu _{1}\left( dy\right) $ (see (\ref{Def_mu})).

As a result, we may apply Lemma \ref{L_Subcr} with our $\delta $ and
conclude that, as $j\rightarrow \infty $%
\begin{equation}
\frac{\mathbf{E}\left[ \vec{\Psi}_{K-x}(\mathbf{W}_{j},\tilde{B}%
_{j}(u),S_{j})e^{S_{j}};\tau _{j}=j\right] }{\mathbf{E}\left[ e^{S_{j}};\tau
_{j}=j\right] }\rightarrow h_{\phi }(u,K-x)  \label{Posit222}
\end{equation}%
for all $v\in \lbrack 0,J]$, where%
\begin{eqnarray*}
&&h_{\phi }(u,K-x) \\
&:&=\iiint \int \vec{\Psi}_{K-x}(\mathbf{w},b,-y)\mathbf{P}_{y}^{+}\left(
\mathbf{W}_{\infty }\in d\mathbf{w}\right) \mathbf{P}^{-}\left( B_{\infty
}(u)\in db\right) \mu _{1}(dy) \\
&=&\iiint \int \phi (w)\Psi _{K-x}(w_{1},b,-y)\mathbf{P}_{y}^{+}\left(
\mathbf{W}_{\infty }\in d\mathbf{w}\right) \mathbf{P}^{-}\left( B_{\infty
}(u)\in db\right) \mu _{1}(dy)
\end{eqnarray*}%
is a bounded function.

Taking $\phi (w)\equiv 1$ and using indentity (\ref{Identit}) and
convergence in (\ref{Posit222}) we conclude that $h(u,K-x)$ in (\ref%
{Posit222aaa}) is given by the relation
\begin{eqnarray}
&&h(u,K-x)  \notag \\
&=&\iiint \Psi _{K-x}(w_{1},b,-y)\mathbf{P}_{y}^{+}\left( W^{+}\in
dw_{1}\right) \mathbf{P}^{-}\left( B_{\infty }(u)\in db\right) \mu _{1}(dy)
\label{Defh}
\end{eqnarray}%
(see also formula (7.12) in \cite{VDD2024}). It was shown in the proof of
Lemma 3.4 in \cite{ABGV2011} that if conditions $B0$ and $B2$ are valid,
then $B_{\infty }(u)<1$ $\mathbf{P}^{-}$-a.s. Thus, $h(u,K-x)>0$.

Further, in view of the estimate
\begin{eqnarray*}
&&\mathbf{E}\left[ \phi \left( \mathcal{Y}^{\theta n,n}\right) \left(
1-F_{(1-\theta )n,j}^{Z_{(1-\theta )n}}(u)\right) ;S_{j}\leq K-x,\tau _{j}=j%
\right] \\
&&\qquad\leq\mathbf{E}\left[ \left( 1-F_{(1-\theta )n,j}^{Z_{(1-\theta
)n}}(u)\right) ;S_{j}\leq K-x,\tau _{j}=j\right] \\
&&\qquad \qquad =\mathbf{E}\left[ 1-F_{0,j}(u);S_{j}\leq K-x,\tau _{j}=j%
\right]
\end{eqnarray*}%
and relations (\ref{Posit222aaa}) and (\ref{Bound_h}) we have
\begin{equation}
h_{\phi }(u,K-x)\leq h(u,K-x)\leq C_{1}(1-u)e^{(K-x)/2}.  \label{Bound-Hphi}
\end{equation}

According to Corollary \ref{C_h}
\begin{eqnarray*}
&&\mathbf{E}\left[ \phi \left( \mathcal{Y}^{\theta n,j+v-\theta n}\right)
\left( 1-F_{j+v-\theta n,j}^{Z_{j+v-\theta n}}(u)\right) ;S_{j}\leq K-x,\tau
_{j}=j\right]  \\
&&\,\leq C\mathbf{E}\left[ 1-F_{0,j}(u);S_{j}\leq K-x,\tau _{j}=j\right]
\leq C(1-u)b_{j}e^{(K-x)/2}.
\end{eqnarray*}%
Hence, applying the dominated convergence theorem we conclude that
\begin{eqnarray}
&&\lim_{j\rightarrow \infty }\int_{0}^{\infty }\int_{0}^{1}\mathbf{P}\left(
F_{0,v}(0)\in du,S_{v}\in dx,L_{v}\geq 0\right)   \notag \\
&&\quad \times \frac{\mathbf{E}\left[ \phi \left( \mathcal{Y}^{\theta
(j+v),j+v}\right) \left( 1-F_{j+v-\theta n,j}^{Z_{(j+v)-\theta n}}(u)\right)
;S_{j}\leq K-x,\tau _{j}=j\right] }{\mathbf{E}\left[ e^{S_{j}};\tau _{j}=j%
\right] }  \notag \\
&&\qquad =\int_{0}^{\infty }\int_{0}^{1}\mathbf{P}\left( F_{0,v}(0)\in
du,S_{v}\in dx,L_{v}\geq 0\right) h_{\phi }(u,K-x)  \notag \\
&&\qquad \quad =\mathbf{E}\left[ h_{\phi }(F_{0,v}(0),K-S_{v});L_{v}\geq 0%
\right] .  \label{ASingleTermRight}
\end{eqnarray}%
Using this fact and summing (\ref{ASingleTermRight}) over $v$ from $0$ to $J$
we get
\begin{eqnarray*}
&&G_{right}\left( K,\phi \right) :=\lim_{J\rightarrow \infty
}\lim_{n\rightarrow \infty }\frac{\mathbf{E}\left[ \phi \left( \mathcal{Y}%
^{\theta n,n}\right) ;\mathcal{R}(n,K),\tau _{n}\in \lbrack n-J+1,n\right] }{%
b_{n}} \\
&&\quad =g_{\alpha ,\beta }(0)\int_{0}^{\infty }e^{-y}U(y)dy\times  \\
&&\qquad \times \int_{0}^{\infty }\int_{0}^{1}\sum_{v=0}^{\infty }\mathbf{P}%
\left( F_{0,v}(0)\in du,S_{v}\in dx,L_{v}\geq 0\right) h_{\phi }(u,K-x) \\
&&\quad \qquad =\int \phi (w)\lambda _{right}(dw),
\end{eqnarray*}%
where%
\begin{eqnarray*}
\lambda _{right}(dw) &=&g_{\alpha ,\beta }(0)\int_{0}^{\infty }e^{-y}U(y)dy
\\
&&\times \mathbf{E}\left[ \sum_{v=0}^{\infty }h_{\phi
}(F_{0,v}(0),K-S_{v})\lambda _{S_{v},F_{0,v}(0)}^{\ast }(dw)\mathbf{I}%
\left\{ L_{v}\geq 0\right\} \right]
\end{eqnarray*}%
with%
\begin{equation*}
\lambda _{x,u}^{\ast }(dw)=\int \int \Psi _{K-x}(w_{1},b,-y)\mathbf{P}%
_{y}^{+}\left( \mathbf{W}_{\infty }\in d\mathbf{w}\right) \mathbf{P}%
^{-}\left( B_{\infty }(u)\in db\right) \mu _{1}(dy)
\end{equation*}%
and
\begin{equation}
\int \lambda _{right}(dw)=G_{right}(K).  \label{massRight}
\end{equation}%
Finally, $G_{right}(K,\phi )\leq G_{right}(K)<\infty $ in view of (\ref%
{Bound-Hphi}).

Thus,
\begin{equation*}
\lim_{n\rightarrow \infty }\frac{\mathbf{E[}\phi \left( \mathcal{Y}^{\theta
n,n}\right) ;\mathcal{R}(n,K)]}{b_{n}}=G_{left}(K,\phi )+G_{right}(K,\phi ),
\end{equation*}%
and%
\begin{eqnarray*}
\lim_{n\rightarrow \infty }\mathbf{E[}\phi \left( \mathcal{Y}^{\theta
n,n}\right) |\mathcal{R}(n,K)] &=&\frac{G_{left}\left( K,\phi \right)
+G_{right}\left( K,\phi \right) }{G_{left}(K)+G_{right}(K)} \\
&=&\int \phi (w)\lambda (dw),
\end{eqnarray*}%
where $\lambda $ is a measure on $D[0,1]$ given by%
\begin{equation*}
\lambda (dw)=\frac{\lambda _{left}(dw)+\lambda _{right}(dw)}{%
G_{left}(K)+G_{right}(K)}.
\end{equation*}%
It is not difficult to check that $\lambda (dw)$ puts entire mass on
strictly positive constant functions and it is a probability measure in view
of (\ref{massLeft}) and (\ref{massRight}).

Theorem \ref{T_generalcase} is proved.


\begin{thebibliography}{99}
\bibitem{agkv} V.I. Afanasyev, J. Geiger, G. Kersting, V. A. Vatutin,
Criticality for branching processes in random environment, \textit{Ann.
Probab.} \textbf{33}:2 (2005), 645--673.

\bibitem{ABGV2011} V.I. Afanasyev, Ch. B\"{o}inghoff, G. Kersting,
V.\thinspace A. Vatutin, Limit theorems for weakly subcritical branching
processes in a random environment, \textit{J.Theor. Probab.} \textbf{25}:3,
(2012), 703--732.

\bibitem{AK1971} K. B. Athreya, S. Karlin, On Branching Processes with
Random Environments: I: Extinction Probabilities, Ann. Math. Statist. 42(5):
1499-1520.

\bibitem{Bil99} P.Billingsley, Convergence of Probability Measures. 2nd
Edition, Wiley, 1999. 277 pp.

\bibitem{Do85} R. A. Doney, \textit{\ Z. Wahrsch. Verw. Gebiete}, \textbf{70}
(1985), 351--360.

\bibitem{Du78} R. Durrett, Conditional limit theorems for some null
reccurent Markov processes, \textit{Ann. Prob.}, \textbf{6} (1978), 798--828.

\bibitem{KV2017} G. Kersting, V. Vatutin, \textit{Discrete Time Branching
Processes in Random Environment,} ISTE and John Wiley \& Sons Limited,
London, 2017.

\bibitem{VD2022} V. A. Vatutin, E. E. Dyakonova, Critical branching
processes evolving in a unfavorable random environment, \textit{Discrete
Math. Appl.}, \textbf{34}:3 (2024), 175--186

\bibitem{VD2023} V.A. Vatutin, E.E. Dyakonova, Population size of a critical
branching processes evolving in an unfavorable environment, \textit{Theory
Probab. Appl.}, \textbf{68}:3 (2023), 411--430.

\bibitem{VDD2023} V.A. Vatutin, K. Dong, E.E. Dyakonova, Random walks
conditioned to stay nonnegative and branching processes in an unfavourable
environment, \textit{Sb. Math.}, \textbf{214}:11 (2023), 1501--1533.

\bibitem{VW2010} V.A. Vatutin, V. Wachtel, Local probabilities for random
walks conditioned to stay positive, \textit{Probab. Theory Related Fields},
\textbf{143}:1-2 (2009), 177--217.

\bibitem{VDD2024} V. A. Vatutin, C. Dong, E. E. Dyakonova, Some functionals
for random walks and critical branching processes in an extremely
unfavourable random environment, \textit{Sb. Math.}, \textbf{215}:10 (2024),
1321--1350.
\end{thebibliography}
\end{document}